\documentclass{article}

\usepackage{arxiv}

\usepackage[utf8]{inputenc} 
\usepackage[T1]{fontenc}    
\usepackage{hyperref}       
\usepackage{url}            
\usepackage{booktabs}       
\usepackage{amsfonts}       
\usepackage{nicefrac}       
\usepackage{microtype}      
\usepackage{lipsum}

\usepackage{graphicx}
\usepackage{multirow}
\usepackage{subcaption}
\usepackage{amsfonts}
\usepackage{amsmath}
\usepackage{amssymb,amsthm,amscd,amstext}
\usepackage[ruled,linesnumbered]{algorithm2e}
\SetKwRepeat{Do}{do}{while}
\usepackage{setspace}
\usepackage[autolanguage]{numprint}
\usepackage[table]{xcolor}
\usepackage{boldline}

\usepackage{framed} 
\usepackage{multicol}
\usepackage{nomencl}
\makenomenclature
\usepackage{mathrsfs} 
\usepackage[numbers]{natbib}

\title{Global optimization using mixed integer quadratic programming on non-convex two-way interaction truncated linear multivariate adaptive regression splines}

\author{
  Xinglong Ju\thanks{This is to indicate the corresponding author.}\\
  Department of Industrial, Manufacturing, \&  Systems Engineering\\
  The University of Texas at Arlington\\
  Arlington, TX 76019, USA \\
  \texttt{xinglong.ju@mavs.uta.edu}\\  \And
 Jay M.\ Rosenberger\\
  Department of Industrial, Manufacturing, \&  Systems Engineering\\
  The University of Texas at Arlington\\
  Arlington, TX 76019, USA\\
  \texttt{xinglong.ju@mavs.uta.edu}\\  \And
    Victoria C.\ P.\ Chen\\
  Department of Industrial, Manufacturing, \&  Systems Engineering\\
  The University of Texas at Arlington\\
  Arlington, TX 76019, USA\\
  \texttt{vchen@uta.edu}\\      \And
Feng Liu\\
  Department of Anesthesia, Critical Care and Pain Medicine\\
  Massachusetts General Hospital, Harvard Medical School\\
  Boston, MA 02114, USA \\
  The Picower The Picower Institute for Learning and Memory\\
	Massachusetts Institute of Technology\\
	Cambridge, MA 02139, USA\\
  \texttt{fliu0@mgh.harvard.edu,\ fengliu@mit.edu}\\
}

\begin{document}
\maketitle

\begin{abstract}
Multivariate adaptive regression splines (MARS) is a flexible statistical modeling
method that has been popular for data mining applications. MARS has also been 
employed to approxmiate unknown relationships in optimzation for complex systems, 
including surrogate optimization, dynamic programming, and two-stage stochastic 
programming.  Given the increasing desire to optimize real world systems, this 
paper presents an approach to globally optimize a MARS model that allows up to 
two-way interaction terms that are products of truncated linear univariate functions 
(TITL-MARS).  Specifally, such a MARS model consists of linear and quadratic structure. 
This structure is exploited to formulate a mixed integer quadratic programming problem 
(TITL-MARS-OPT).
To appreciate the contribution of TITL-MARS-OPT, one must recognize that popular
heurstic optimization approaches, such as evolutionary algorithms, do not guarantee
global optimality and can be computationally slow. The use of MARS maintains the
flexibility of modeling within TITL-MARS-OPT while also taking advantage of
the linear modeling structure of MARS to enable global optimality.
Computational results compare TITL-MARS-OPT with a genetic algorithm for two types of cases.
First, a wind farm power distribution case study is described and then other TITL-MARS forms 
are tested. The results show the superiority of TITL-MARS-OPT over the genetic algorithm 
in both accuracy and computational time.
\end{abstract}


\section{Introduction}
\begin{figure}[h]
	\begin{center}
		\includegraphics[width=2.5in]{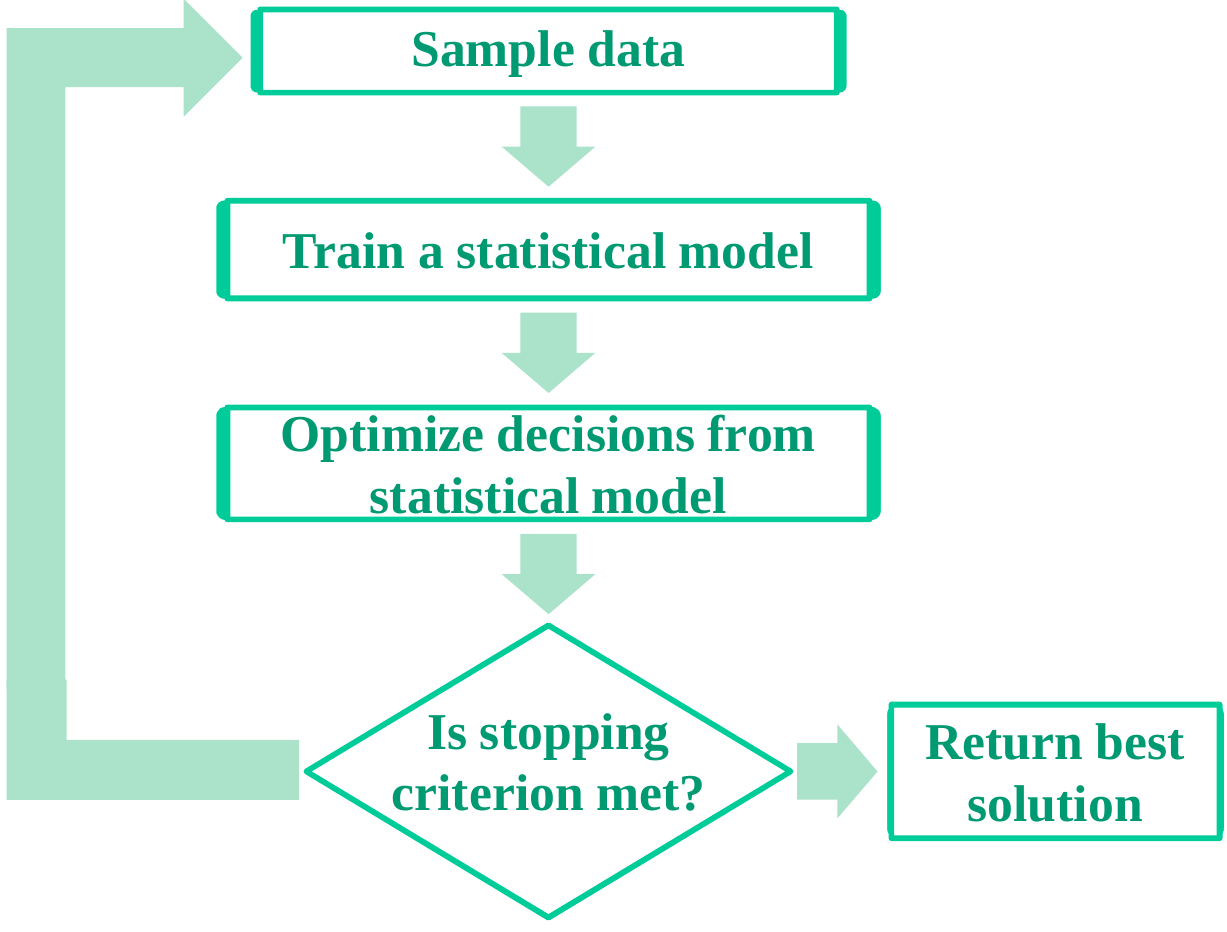}
		\caption{Optimization of complex systems}
		\label{fig:ocs_fc}
	\end{center}
\end{figure}

Optimization for complex systems often involves fitting a system prediction model to estimate how a system performs and then optimizing the decisions based on the system prediction model as shown in Figure~\ref{fig:ocs_fc}. Two major tasks in optimization of complex systems include training or meta-modeling a statistical or system model and optimizing input or decisions based on statistical model.

In real world complex systems, underlying relationships are commonly unknown and are approximated 
from data using empirical models.  \citet{wu2004travel} applied support vector regression in travel time prediction and proved support vector regression was applicable in traffic data analysis. \cite{lv2014traffic} applied a deep learning approach with autoencoders in traffic flow prediction. \cite{ohlmacher2003using} used logistic regression to generate a landslide-hazard map to predict landslide hazards. \cite{leathwick2005using} used multivariate adaptive regression splines to predict the distributions of freshwater diadromous fish.

%
%

If one seeks to optimize a complex system, the optimization
method would need to be able to handle the data-driven approximation models.  Given the wide range
of possible approximation models, such as machine learning algorithms, the most commonly
employed optimization approach in these situations is a heuristic approach, such as an evolutionary
algorithm, that cannot guarantee global optimality.  Rather than having the approximation model
dictate the need for a heuristic optimization method, the research in this paper seeks a balance that
utilizes a flexible approximation model with structure that can be exploited to enable 
true global optimization.  In other words, the ``best of both worlds'' is sought, by achieving global
optimality while still maintaining a flexible approximation model.
The approxmation model of choice in this paper is multivariate adaptive regression splines (MARS),
introduced by machine learning pioneer Jerome Friedman in 
1991~\citep{friedman1991multivariate}.  The structure of MARS
is based on a linear statistical spline model and provides a flexible fit to data while also
achieving a parsimonious model. 


The desire to conduct global optimization is seen in many applications, and there are a number
of approaches classified as global optimization methods~\citep{horst2000introduction}.  The primary
challenge in achieving global optimality is that many real world applications involve multiple
local optima.  Finding a global optimum requires sifting through the local optima and recognizing
when one is suboptimal.  The vast majority of applications employ heuristic search algorithms
seek to overcome the challenge of local optima, but do not guarantee global optimality.
Examples include heuristics based on evolutionary algorithms~\citep{morra2018optimization,russo2018knowledge,dahi2018stop},
particle swarm optimization~\citep{alswaitti2018density},
the grasshopper optimization method~\citep{ewees2018improved}, and
the weighted superposition attraction method~\citep{baykasouglu2018dynamic}. 
In order to guarantee global optimality, the approach in this paper takes advantage of well-known
properties of mixed integer and quadratic programming (MIQP)~\citep{bliek1u2014solving}.
%


Some recent applications in which MARS has been employed for empirical modeling include
a water pollution prediction problem~\citep{kisi2016application},
the head load in a building~\citep{roy2018estimating}, 
the estimation of landfill leachate~\citep{bhatt2017estimating},
and the damage identification for web core composite bridges~\citep{mukhopadhyay2018multivariate}, 
In optimization problems, MARS has been employed as the empirical model to approximate
unknown relationships in a variety of applications.  For stochastic dynamic programming,
the use of MARS to approximate the value function was introduced by~\cite{chen1999applying}.
Since then, the MARS value function approximation approach has been used to numerically solve
a 30-dimensional water reservoir management problem~\citep{cervellera2006optimization},
a 20-dimensional wastewater treatment system~\citep{tarun2011optimizing,tsai2004stochastic,tsai2005flexible},
and a 524-dimensional nonstationary ground-level ozone pollution control problem~\citep{yang2009decision}.
%
In revenue management, MARS was employed to estimate upper and lower bounds for the
value function of a Markov decision problem~\citep{chen2003solving,siddappa2007refined},
and MARS was used to represent the revenue function in airline overbooking
optimization~\cite{siddappa2008optimising}.
%
In two-stage stochastic programming, MARS was used to efficiently represent the
expected profit function for an airline fleet assignment problem~\citep{pilla2008statistical}.
This fleet assignment research was extended to utilize a cutting plane method with MARS
to conduct the optimization~\citep{pilla2012multivariate}.

%
The contribution of this current work extends the approach of 
\cite{martinez2017global},
who developed a piece-wise linear MARS structure and formulated a mixed integer and
linear programming problem to globally optimize vehicle design parameters to
improve performance in crash simulations.
The piece-wise linear MARS function may be nonconvex, and the approach of Martinez et al.
will yield a global optimum.  However, restricting to piece-wise linear forms limits
the flexibility of the empirical model.  Hence, in the current work, the MARS form
employed is based on the original MARS model.  The primary challenge for an optimization
method is handling the nonconvex MARS interaction terms, which are products of univariate terms.
By restricting to two-way interactions, we can utilize quadratic programming methods.
In real world applications, two-way interactions are commonly sufficient for empirical modeling~\cite{kutner2005applied}

In summary, the contribution of the presented approach is a MIQP
global optimization method for a MARS model that allows up to 
two-way interaction terms that are products of truncated linear univariate functions 
(TITL-MARS).  This approach is referred to as TITL-MARS-OPT and is compared against
a genetic algorithm for two types of cases.
First, a wind farm power distribution case study is described, and then other TITL-MARS forms 
are tested.
Python code for TITL-MARS and TITL-MARS-OPT will be made available on GitHub upon 
acceptance of this paper~(\url{https://github.com/JuXinglong/TITL-MARS-OPT}). 

The rest of this paper is organized as follows. 
Section~\ref{background} describes background on TITL-MARS. 
Section~\ref{formulation} presents the MIQP formulation for TITL-MARS-OPT. 
The computational study is given in Section~\ref{experiments}, and 
Section~\ref{conclusion} concludes the paper.

\section{Background of two-way interaction truncated linear multivariate adaptive regression splines}\label{background}


This section introduces the two-way interaction truncated linear MARS (TITL-MARS) model.
The two-way interaction truncated linear MARS regression model with the response variable ${f}( \mathbf{x}_i)$ is to be built on the independent variable $\mathbf{x}_i$ and can be written in the form of the linear combination of the basis functions as~\cite{friedman1991multivariate}
\begin{equation}\label{eq:MARS_2way}
\begin{array}{l}
\hat{f}( \mathbf{x})=a_0+\sum_{m=1}^{M}\left\{a_m \cdot B_m(\mathbf{x}) \right\}.
\end{array}
\end{equation}
The MARS model is denoted as $\hat{f}( \mathbf{x})$, and $a_0$ is the constant term of the model. The basis function is denoted as $B_m(\mathbf{x})$, and $a_m$ is the coefficient of $B_m(\mathbf{x})$. The index of the basis function is denoted as $m$, and $M$ is the total number of basis functions. The basis function $B_m(\mathbf{x})$ using the truncated linear term has the following form
\begin{equation}\label{eq:MARS_2way_B}
\begin{array}{l}
B_m(\mathbf{x})= \prod_{k=1}^{Km}[s_{k,m} \cdot (x_{v(k,m)}-t_{v(k,m)}) ]_+.
\end{array}
\end{equation}
The truncated linear term is denoted as ${[s_{k,m} \cdot (x_{v(k,m)}-t_{v(k,m)}) ]}_+$, and the basis function $B_m(\mathbf{x})$ is the product of truncated linear terms.
The index of the truncated linear term in $B_m(\mathbf{x})$ is denoted as $k$, and $K_m$ is the total number of truncated linear terms in $B_m(\mathbf{x})$. The sign of the truncated linear term is $s_{k,m}$, which can be $+1$ or $-1$. The $v$-th component of $\mathbf{x}$ is denoted as $x_{v(k,m)}$, and $t_{v(k,m)}$ is the corresponding knot value.
\textit{TITL-MARS is the special case of MARS in which $K_m \leqslant 2$.}

\section{Formulation of two-way interaction truncated linear MARS using mixed integer quadratic programming}\label{formulation}

The general mixed integer quadratic programming problem~\citep{bliek1u2014solving} is given as
\begin{equation}
\label{eq:miqpp}
\begin{array}{ll}
\min & \frac{1}{2} \textbf{z}^T\textbf{Q}\textbf{z}+\textbf{c}^T\textbf{z} \\
{s.t.} & \textbf{A}\textbf{z}=\textbf{b} \\
{} & \textbf{l} \leqslant \textbf{z} \leqslant \textbf{u}\\
{} &\textbf{z} \in \mathbb{R}^P \times \mathbb{Z}^{D-P},\\
\end{array}
\end{equation}
while the two-way interaction truncated linear MARS is given in Section~\ref{background}. In~\eqref{eq:miqpp}, the decision variable is $\textbf{z}$, and the quadratic coefficients matrix is $\textbf{Q}$. The coefficients of the linear terms in the objective funtion are in vector $\textbf{c}$.  The linear constraints are denoated as $\textbf{A}\textbf{z}=\textbf{b}$. The lower bound and upper bound of $\textbf{z}$ are $\textbf{l}$ and $\textbf{u}$, respectively. The dimension of $\textbf{z}$ is $D$. There are $P$ dimensions of real values, and $D-P$ dimensions of integers.
Problem in the form~\ref{eq:miqpp} can be solved using the CPLEX solver. 

The TITL-MARS optimization prolbem is given as follows. \begin{eqnarray}
\label{eq:mars_opt_1}
\min\,\, & \hat{f}( \mathbf{x})=a_0+\sum_{m=1}^{M}\left\{a_m \cdot B_m(\mathbf{x}) \right\}  \\
\label{eq:mars_opt_2}
\mbox{s.t.}\,\, & \textbf{l} \leqslant \textbf{x} \leqslant \textbf{u}\\
\label{eq:mars_opt_3}
{} &\textbf{x} \in \mathbb{R}^P \times \mathbb{Z}^{D-P}\\
\label{eq:mars_opt_4}
&B_m(\mathbf{x})= \prod_{k=1}^{Km}[s_{k,m} \cdot (x_{v(k,m)}-t_{v(k,m)}) ]_+\\
\label{eq:mars_opt_5}
&{{[s_{k,m} \cdot (x_{v(k,m)}-t_{v(k,m)}) ]}_+}= \max \{s_{k,m} \cdot (x_{v(k,m)}-t_{v(k,m)}),0\}
\end{eqnarray}
The objective function~\eqref{eq:mars_opt_1} is the TITL-MARS model. The constraint set~\eqref{eq:mars_opt_2} is the boundary of $\mathbf{x}$. The constraint set~\eqref{eq:mars_opt_3} specifies the data types. Constraints~\eqref{eq:mars_opt_4} and \eqref{eq:mars_opt_5} specify the basis functions and the truncated linear terms.

Let $\mathcal{M}$ denote an upper bound of $|x_{v(k,m)}-t_{v(k,m)}|$ and $|t_{v(k,m)}-x_{v(k,m)}|$. Let $y_{k,m}$ be an indicator variable for the nonnegativity of $s_{k,m} \cdot (x_{v(k,m)}-t_{v(k,m)})$, an let $\eta_{k,m}$ denote the univariate truncated linear function, given as
\begin{equation}\label{eq:utlbf1}
\eta_{k,m}=[s_{k,m} \cdot (x_{v(k,m)}-t_{v(k,m)})]_+=\max\{s_{k,m} \cdot (x_{v(k,m)}-t_{v(k,m)}),0\}.
\end{equation}
Specifically, when $s_{k,m} \cdot (x_{v(k,m)}-t_{v(k,m)}) \geqslant 0$, $y_{k,m}=1$ and $\eta_{k,m}=s_{k,m} \cdot (x_{v(k,m)}-t_{v(k,m)})$, otherwise $y_{k,m}=0$ and $\eta_{k,m}=0$.

The TITL-MARS optimization prolbem can be formulated into a general mixed integer quadratic programming problem as follows. 
\begin{eqnarray}
\label{eq:mars_opt_formu_1}
\min\,\, & a_0+\sum_{m=1}^{M}\left\{a_m \cdot \prod_{k=1}^{Km} \eta_{k,m} \right\}   \\
\mbox{s.t.}\,\, & s_{k,m} \cdot (x_{v(k,m)}-t_{v(k,m)})\leqslant \eta_{k,m} \leqslant s_{k,m} \cdot (x_{v(k,m)}-t_{v(k,m)})+\mathcal{M}\cdot(1-y_{k,m}),\nonumber \\
\label{eq:mars_opt_formu_3}
{}&\forall k=1,\ldots, K_m, \forall m=1,\ldots, M\\
\label{eq:mars_opt_formu_4}
{}&0\leqslant \eta_{k,m}\leqslant \mathcal{M}\cdot y_{k,m},\forall k=1,\ldots, K_m, \forall m=1,\ldots, M\\
\label{eq:mars_opt_formu_5}
{}&\textbf{l} \leqslant \textbf{x} \leqslant \textbf{u}\\
\label{eq:mars_opt_formu_6}
{} &\textbf{x} \in \mathbb{R}^P \times \mathbb{Z}^{D-P}\\
\label{eq:mars_opt_formu_7}
{}&\eta_{k,m} \in \mathbb{R},\forall k=1,\ldots, K_m, \forall m=1,\ldots, M\\
\label{eq:mars_opt_formu_8}
{}&y_{k,m} \in \mathbb{B},\forall k=1,\ldots, K_m, \forall m=1,\ldots, M.
\end{eqnarray}
The objective~\eqref{eq:mars_opt_formu_1} is the TITL-MARS model. Equations~\eqref{eq:mars_opt_formu_3} -~\eqref{eq:mars_opt_formu_8}  formulate the basis functions into linear constraints and specifying the boundaries and data types.

\begin{figure}[h]
	\begin{center}
		\includegraphics[width=1.8in]{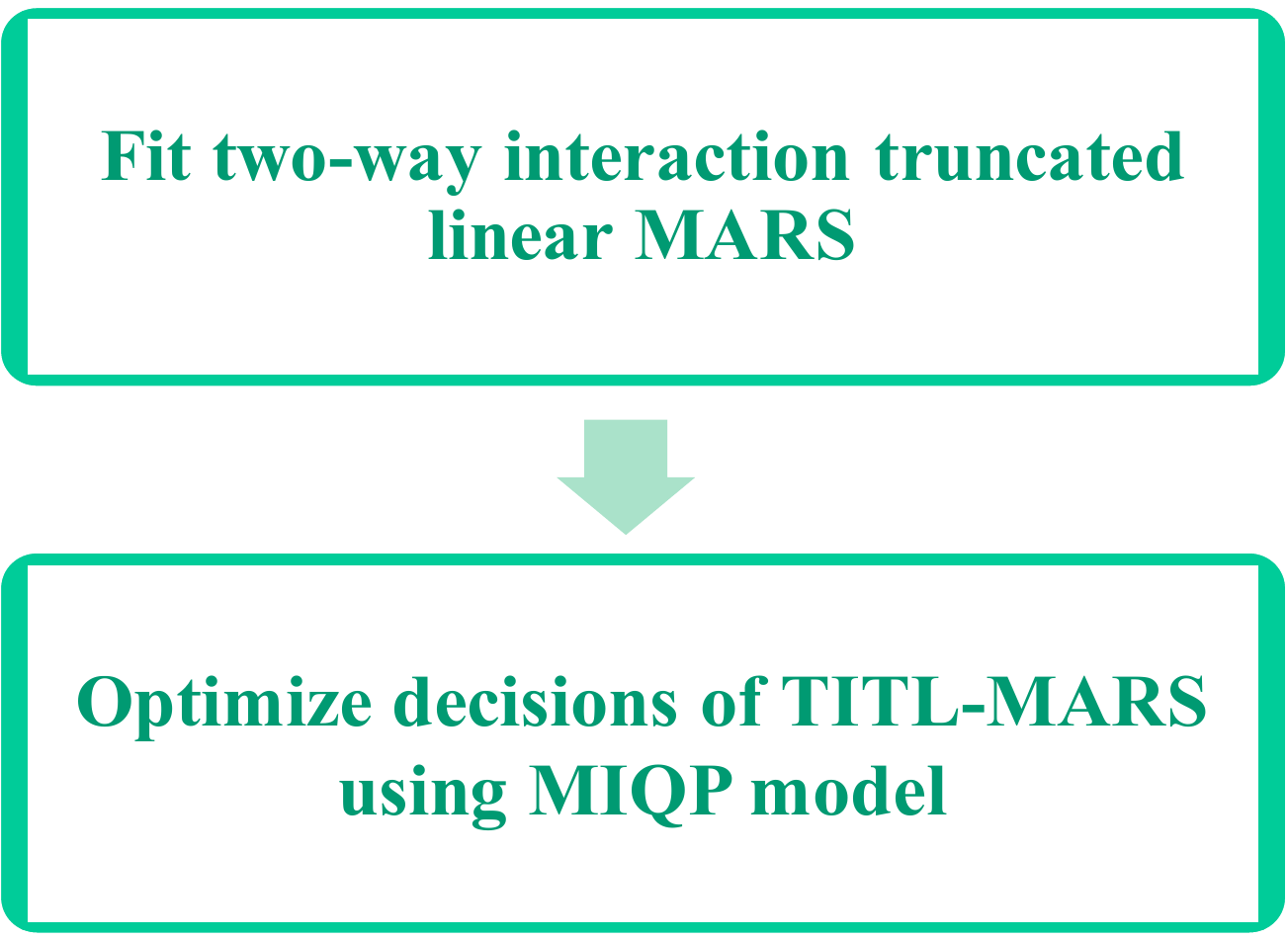}
		\caption{TITL-MARS-OPT optimization process}
		\label{fig:titl_mars_opt_fc}
	\end{center}
\end{figure}

TITL-MARS-OPT is an optimization process as shown in Figure~\ref{fig:titl_mars_opt_fc}. The process has two steps. The first step is to fit TITL-MARS model, and the second step optimizes TITL-MARS model using MIQP. The benefits of the optimization process has two aspects. First, TITL-MARS can be fit using most commercial MARS software. Second, MIQP can be globally optimized using CPLEX~\citep{ilog12ibm}.

\section{Experiments and results}\label{experiments}
In this section, first the genetic algorithm for TITL-MARS optimization is given, and then the presented TITL-MARS-OPT is tested on wind farm power distribution TITL-MARS models and other mathematical models with the genetic algorithm (TITL-MARS-GA) as a benchmark.
\subsection{Genetic algorithm}
The genetic algorithm can also be used as an optimization method to optimize the function (TITL-MARS-GA), as given in Algorithm~\ref{alg:genetic}~\citep{michalewicz1996evolution}, where the input is the two-way interaction MARS model and the maximum generation number $M_{\max}$, and the output is an optimum and an optimum value. In the ``initialization'' step (line 1 in Algorithm~\ref{alg:genetic}), we generate a population and code the individuals from the decimal form to the binary form. In the ``fitness value'' step (line 2), we decode the individuals from the binary form to the decimal form and evaluate each of the individual's decimal values in the MARS model function to obtain the fitness value. In the ``keep the best'' step (line 3), we sort the individuals by their fitness values and store the individual with the best fitness value. The ``selection'' step (line 6) selects parents from the prior population. The ``crossover'' (line 7) chooses two parents and produces a new population. The ``mutation'' (line 8) chooses one point within an individual and changes it from 1 to 0 or from 0 to 1.  In this paper, the TITL-MARS-GA algorithm is used as a benchmark compared with the TITL-MARS-OPT method.
\begin{algorithm}
	\small
	\caption{Genetic algorithm for TITL-MARS optimization}\label{alg:genetic}
	\KwData{$\hat{f}( \mathbf{x})=a_0+\sum_{m=1}^{M}\left\{a_m \cdot \prod_{k=1}^{Km}[s_{k,m} \cdot (x_{v(k,m)}-t_{v(k,m)}) ]_+ \right\}, M_{\max}$}
	\KwResult{$\mathbf{x}_{\max},f(\mathbf{x}_{\max})$ }
	\textbf{Initialization}: Generate a population and code the individuals from decimal to binary.\\
	\textbf{Fitness value}: Decode individuals from binary to decimal and get function value.\\
	\textbf{Keep the best}: Store the individual with highest or lowest fitness value.\\
	$\mbox{gen}=1$\\
	\While{\mbox{gen}$<M_{\max}$}
	{
		\textbf{Selection}: Select parents from prior population.\\
		\textbf{Crossover}: Choose two parents and produce a new population.\\
		\textbf{Mutation}: Choose one point and $1\rightarrow0$ or $0\rightarrow1$.\\
		\textbf{Fitness value}: Decode individuals from binary to decimal and get function value.\\
		\textbf{Keep the best}: Store the individual with highest or lowest fitness value.\\
		$\mbox{gen}=\mbox{gen}+1$
	}
\end{algorithm}

The parameters of TITL-MARS-GA in this paper are from the literatures \cite{grefenstette1986optimization} and \cite{michalewicz1996evolution}, and given in Table~\ref{tab:titi_mars_ga_parameters}.
\begin{table}[h]
	\centering
	\small
	\caption{Parameter settings of TITL-MARS-GA}
	\label{tab:titi_mars_ga_parameters}
	\scalebox{1.0}{
	\begin{tabular}{c|c|c}
		\hline
		Parameter&\cite{grefenstette1986optimization}& \cite{michalewicz1996evolution}\\
		\hline
		Population size&30&50 \\
		Maximum number of generations&	300&1000 \\
		Crossover rate&0.9 & 0.8\\
		Mutation rate&	0.01& 0.15\\
		\hline
	\end{tabular}}
\end{table}

\begin{figure}[h]
	\begin{center}
		\includegraphics[width=1.8in]{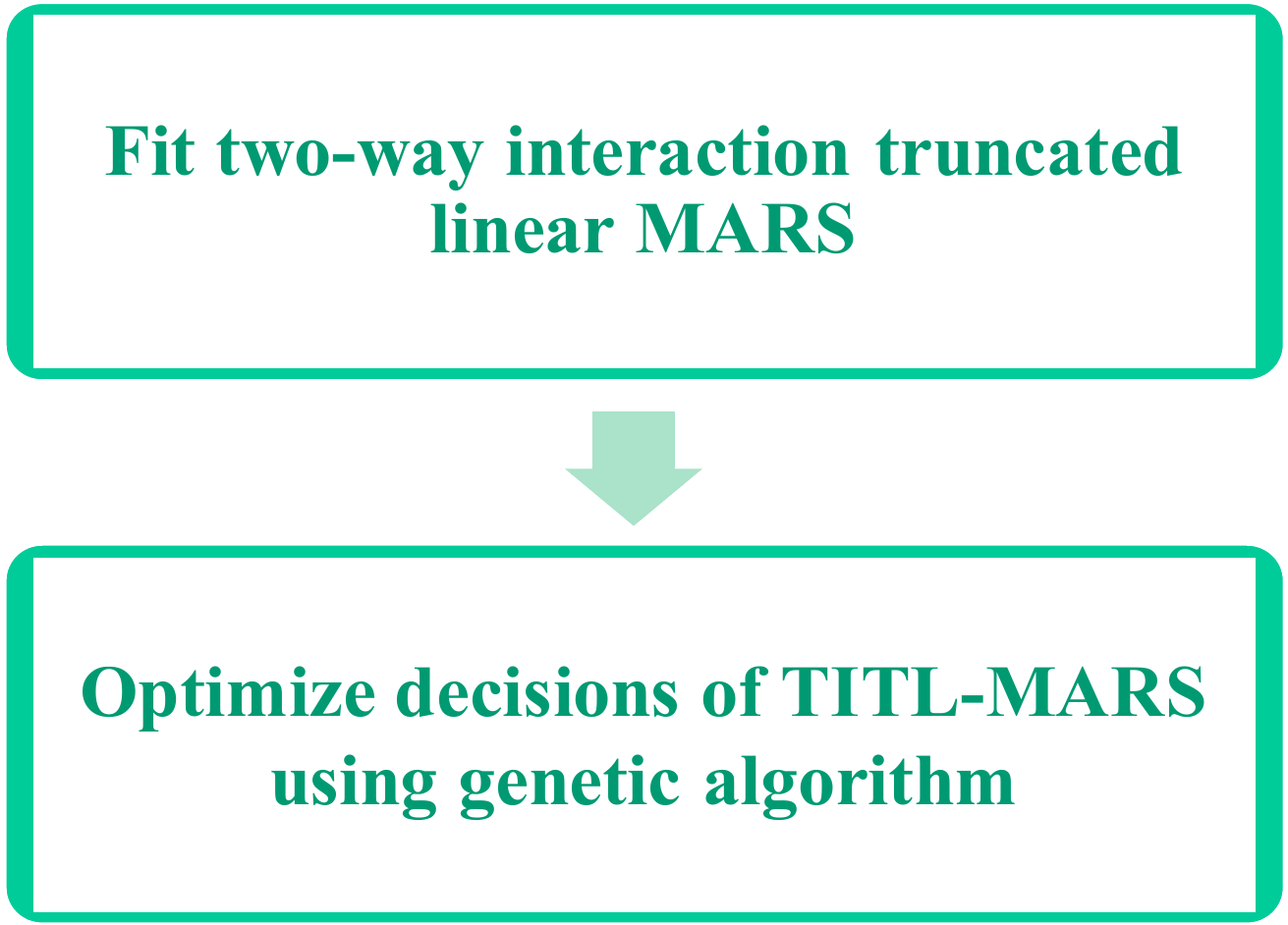}
		\caption{TITL-MARS-GA optimization process}
		\label{fig:titl_mars_ga_fc}
	\end{center}
\end{figure}

The TITL-MARS-GA optimization process has two steps as shown in Figure~\ref{fig:titl_mars_ga_fc}. The first step fits a two-way interaction truncated linear MARS model, and the second step optimizes decisions of TITL-MARS using the genetic algorithm. The drawback of TITL-MARS-GA is that it does not guarantee global optimality.

\subsection{Experimental environment}
The experiments are run on a workstation with 64 bit Windows 10 Enterprise system. The CPU version is an Intel(R) CPU E3-1285 v6 @ 4.10GHz, and the RAM has 32 GB.The programming code is written in Python version is 3.6, and the CPLEX solver version is 12.8.

\subsection{Optimization of wind farm power distribution function}
\begin{figure}[h]
	\begin{center}
		\includegraphics[width=1.8in]{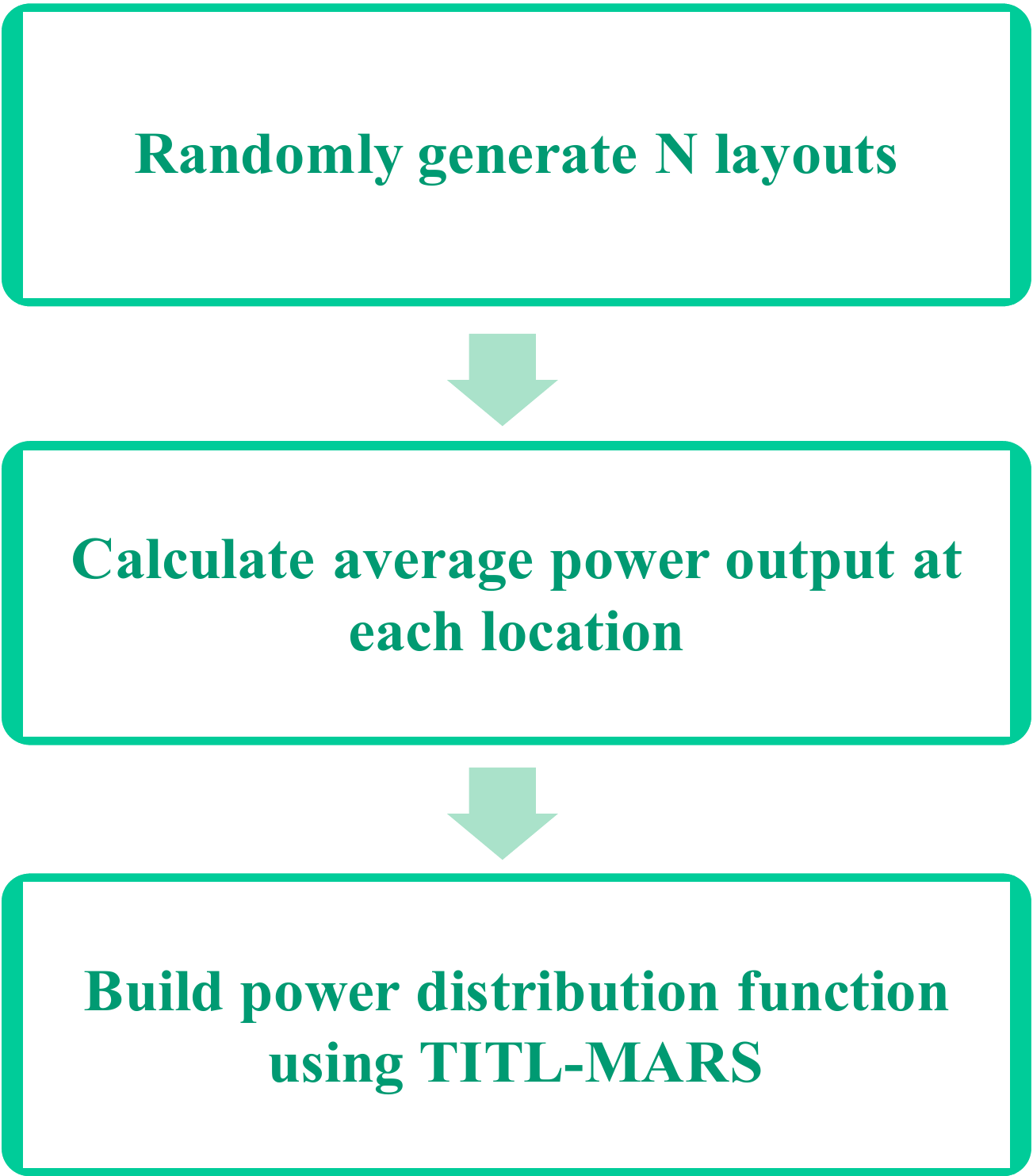}
		\caption{Steps to generate a TITL-MARS wind farm power distribution model}
		\label{fig:titl_mars_wf_fc}
	\end{center}
\end{figure}

Wind farm power is of paramount significance as a renewable energy source. In this paper, the Monte Carlo method~\citep{metropolis1949monte} is used to generate random wind farm layouts, and the  TITL-MARS method is used to study the power distribution under certain wind speeds and directions. After the TITL-MARS model is generated, the TITL-MARS-OPT method is used to study the best turbine position and the worst position. We use the following steps to generate the wind farm power distribution function, as shown in Figure~\ref{fig:titl_mars_wf_fc}. First, we randomly generate $N$ wind farm layouts. Second, we calculate average power output at each location. Third, we use the data from second step to build the TITL-MARS power distribution model.

\begin{figure*}[ht]
	\centering
	\begin{subfigure}[b]{0.4\textwidth}
		\centering
		\includegraphics[height=1.6in]{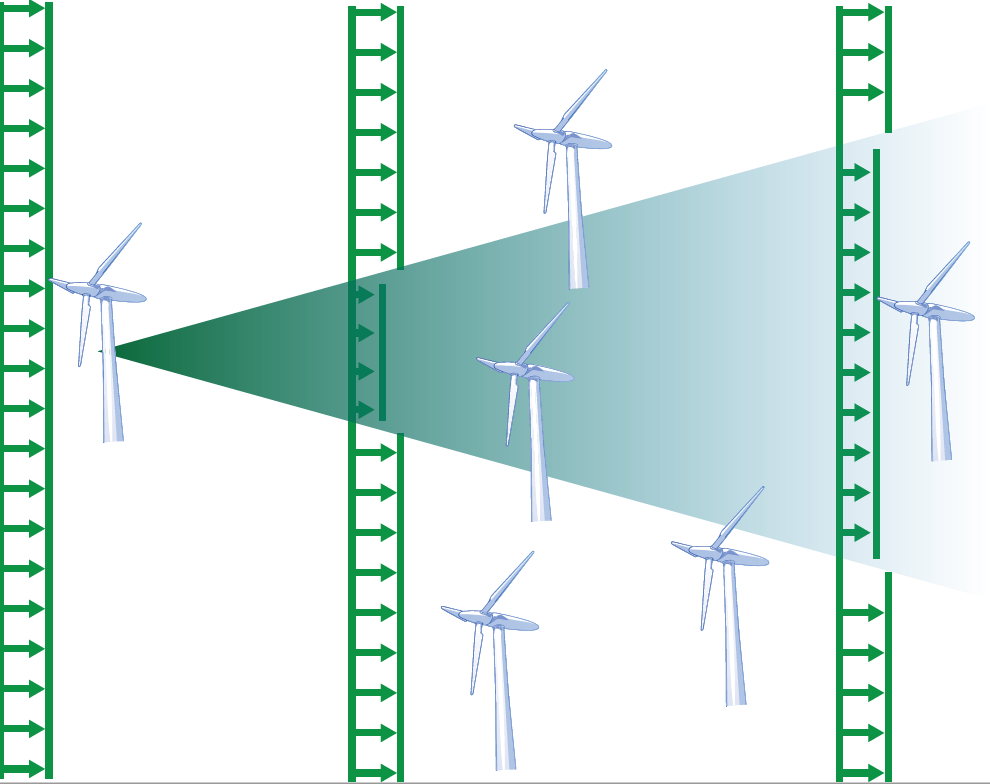}
		\caption[]%
		{{\small Wind wake effect model}}    
		\label{fig:wf_wake_m_1}
	\end{subfigure}
	\hspace{0in}
	\begin{subfigure}[b]{0.5\textwidth}  
		\centering 
		\includegraphics[height=1.6in]{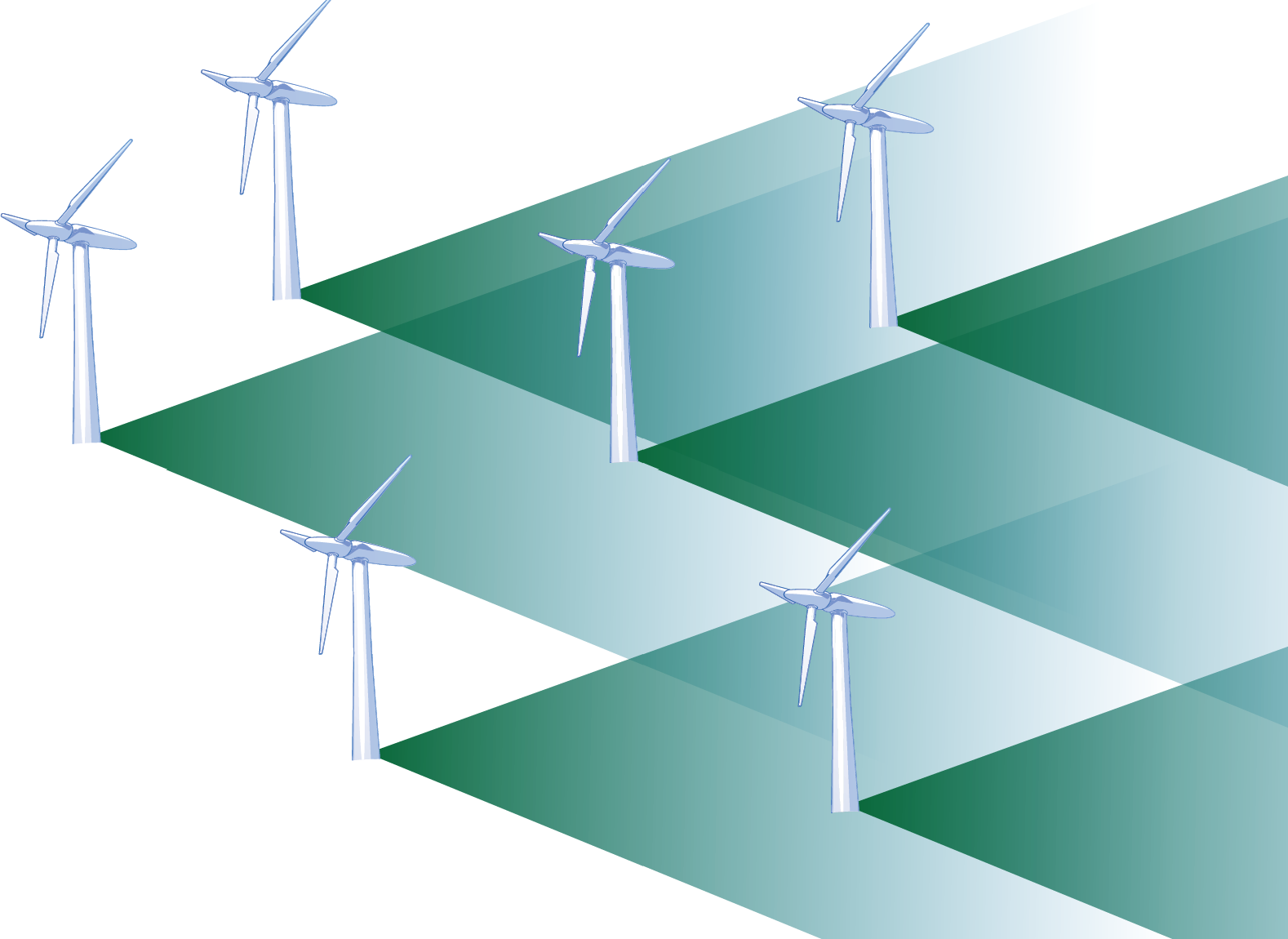}
		\caption{{\small Wake effect illustration}}
		\label{fig:wf_wake_m_2}
	\end{subfigure}
	\caption[]
	{\small Wake effect} 
	\label{fig:wf_wake_m}
\end{figure*}

\begin{figure}[]
	\begin{center}
		\includegraphics[width=3in]{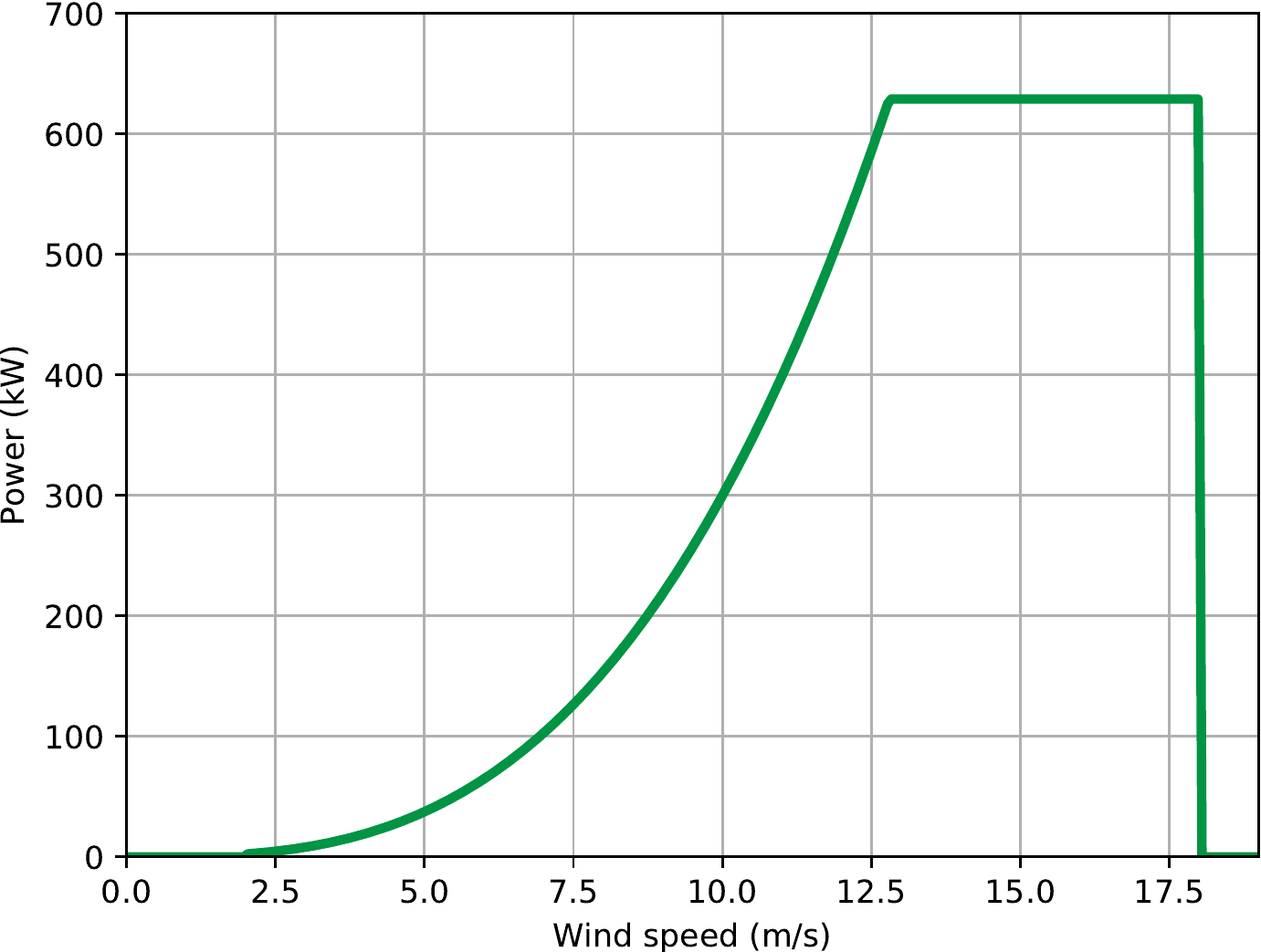}
		\caption{Relationship of wind speed with output power of a wind turbine}
		\label{fig:power_curve}
	\end{center}
\end{figure}

After the wind passes through a wind turbine $j$, a part of the wind energy will be absorbed by turbine $j$ and leave the downstream wind with the reduced speed, which is called the wake effect~\citep{jensen1983note}, and the wake effect model is shown in Figure~\ref{fig:wf_wake_m}. 
Wind speed at turbine $i$ with the wake effect of turbine $j$ is $v_{i,j}$ and can be calculated as
\begin{equation}
\label{eq:rsi}
v_{i,j}=v_0 \left( 1-\frac{2}{3} \cdot \frac{R_j^2}{r_j^2}  \right).
\end{equation}
$R_j$ is the radius of the wind turbine $j$, and $r_j$ is the wake radius of the wind turbine $j$. The final wind speed $v_i$ at turbine $i$ with multiple wake effects is given as
\begin{equation}
\label{eq:vi}
\begin{array}{l}
v_i=v_0 \left[ 1-\sqrt{\sum \limits_{j \in \Phi_{i}} \left(1-\frac{v_{i,j}}{v_0}\right)^2} \right],
\end{array}
\end{equation}
where $\Phi _{i}$ is the index set of the turbines which are upwind of the turbine $i$. Afterwards, the actual power of turbine $i$ can be obtained as~\cite{liu2013electric} 
\begin{equation}\label{eq:powercurvef}
p(v_i) = \left\{ {\begin{array}{*{20}{ll}}
	{0, } & v_{i} < 2\\
	{0.3{v_i}^3,} & {2 \leqslant v_{i} < 12.8{\mkern 1mu} }\\
	{629.1,}  &  12.8{\mkern 1mu} \leqslant {v_{i} \leqslant 18}\\
	{0, } & v_{i} > 18{\mkern 1mu},\\
	\end{array}} \right.
\end{equation}
and the power curve is shown in Figure~\ref{fig:power_curve}. which is the relationship between the wind turbine power and the wind speed.

The wind farm power distribution is generated using the Monte Carlo methods for a given wind farm and a specific wind distribution.

$f_{w1}$ is generated from a  wind farm where there is only one wind speed and one direction. The wind farm is divided into 41 by 41 cells, and each cell has a width of 308 m. The wind is from northeast ($\frac{\pi}{4}$) at 15 m/s.

$f_{w2}$  is generated from a  wind farm where the wind has only one wind speed and four directions. The wind farm has the same dimension as that of $f_{w1}$. The wind is from north (0), south ($\pi$), east ($\frac{\pi}{2}$),  and west ($\frac{3\pi}{2}$) at 15 m/s. 

$f_{w3}$ is generated from a  wind farm where the wind has only one wind speed at 15 m/s and six directions, 0, $\frac{\pi}{3}$, $\frac{2\pi}{3}$, $\pi$, $\frac{4\pi}{3}$, and $\frac{5\pi}{3}$.

$f_{w4}$ is generated from a  wind farm where the wind has three wind speeds, 12 m/s, 10 m/s, and 8 m/s, and 12 directions, 0, $\frac{\pi}{6}$, $\frac{\pi}{3}$, $\frac{\pi}{2}$, $\frac{2\pi}{3}$, $\pi$, $\frac{7\pi}{6}$, $\frac{4\pi}{3}$, $\frac{3\pi}{2}$, $\frac{5\pi}{3}$, and $\frac{11\pi}{6}$.

\begin{figure*}[]
	\centering
	\begin{subfigure}[b]{0.45\textwidth}  
		\centering 
		\includegraphics[width=2.2in]{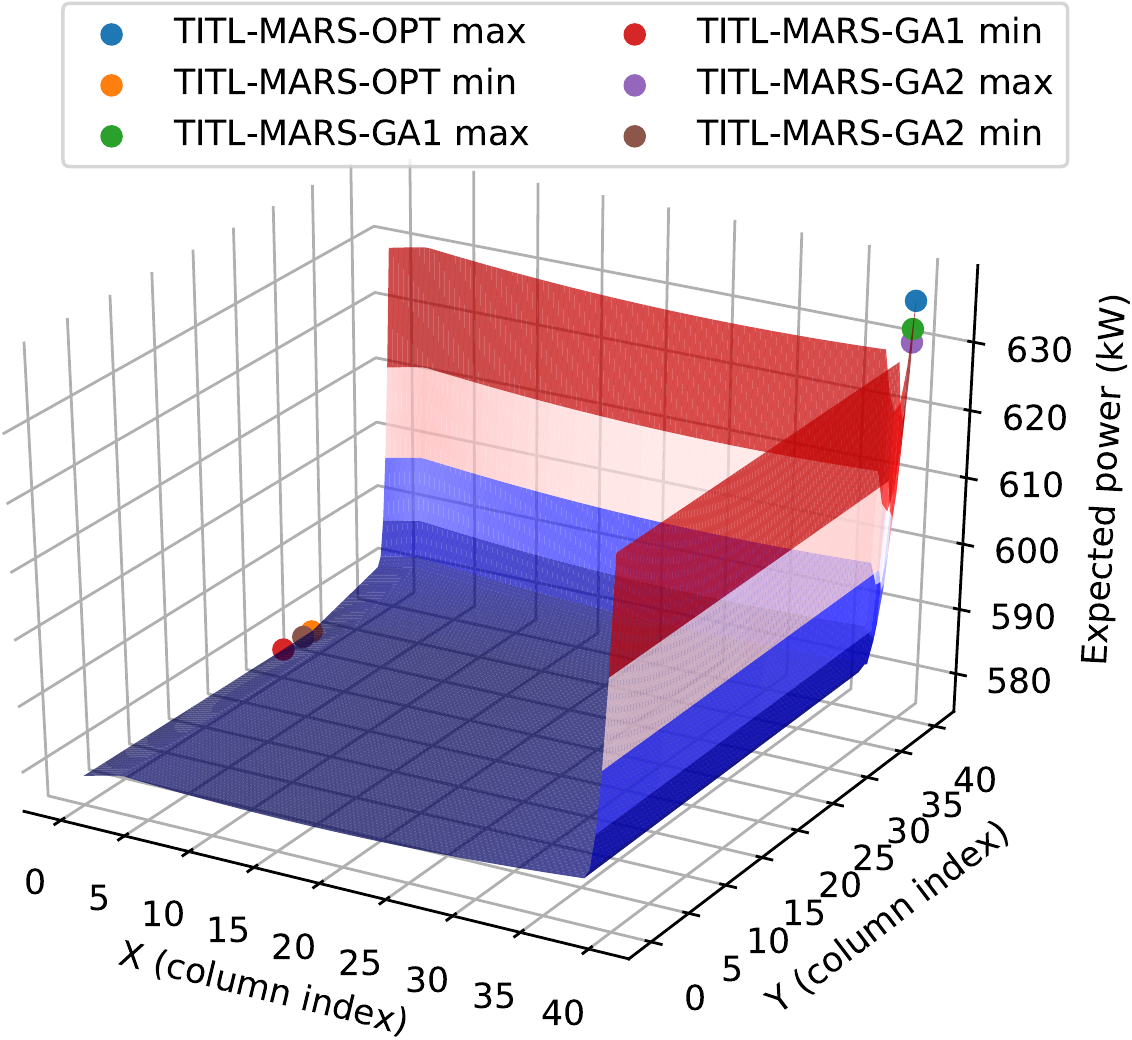}
		\caption[]%
		{{\small Case 1 : $f_{w1}$}}
		\label{fig:mwd_surface_ex_1}
	\end{subfigure}
	\hspace{0em}
	\begin{subfigure}[b]{0.45\textwidth}
		\centering
		\includegraphics[width=2.2in]{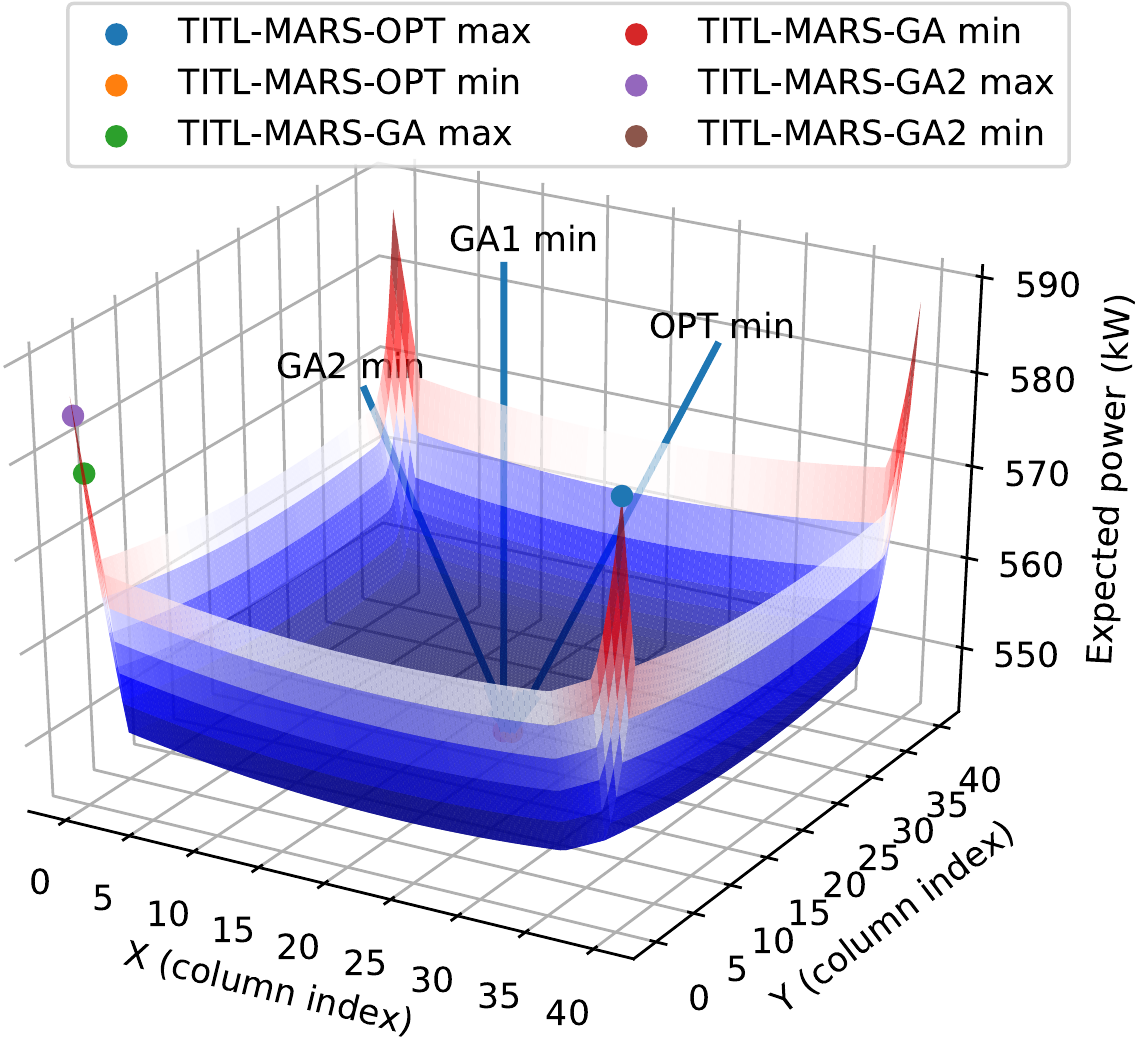}
		\caption[]%
		{{\small Case 2 : $f_{w2}$}}    
		\label{fig:mwd_surface_ex_2}
	\end{subfigure}
	\vspace{0\baselineskip}
	\begin{subfigure}[b]{0.45\textwidth}  
		\centering 
		\includegraphics[width=2.2in]{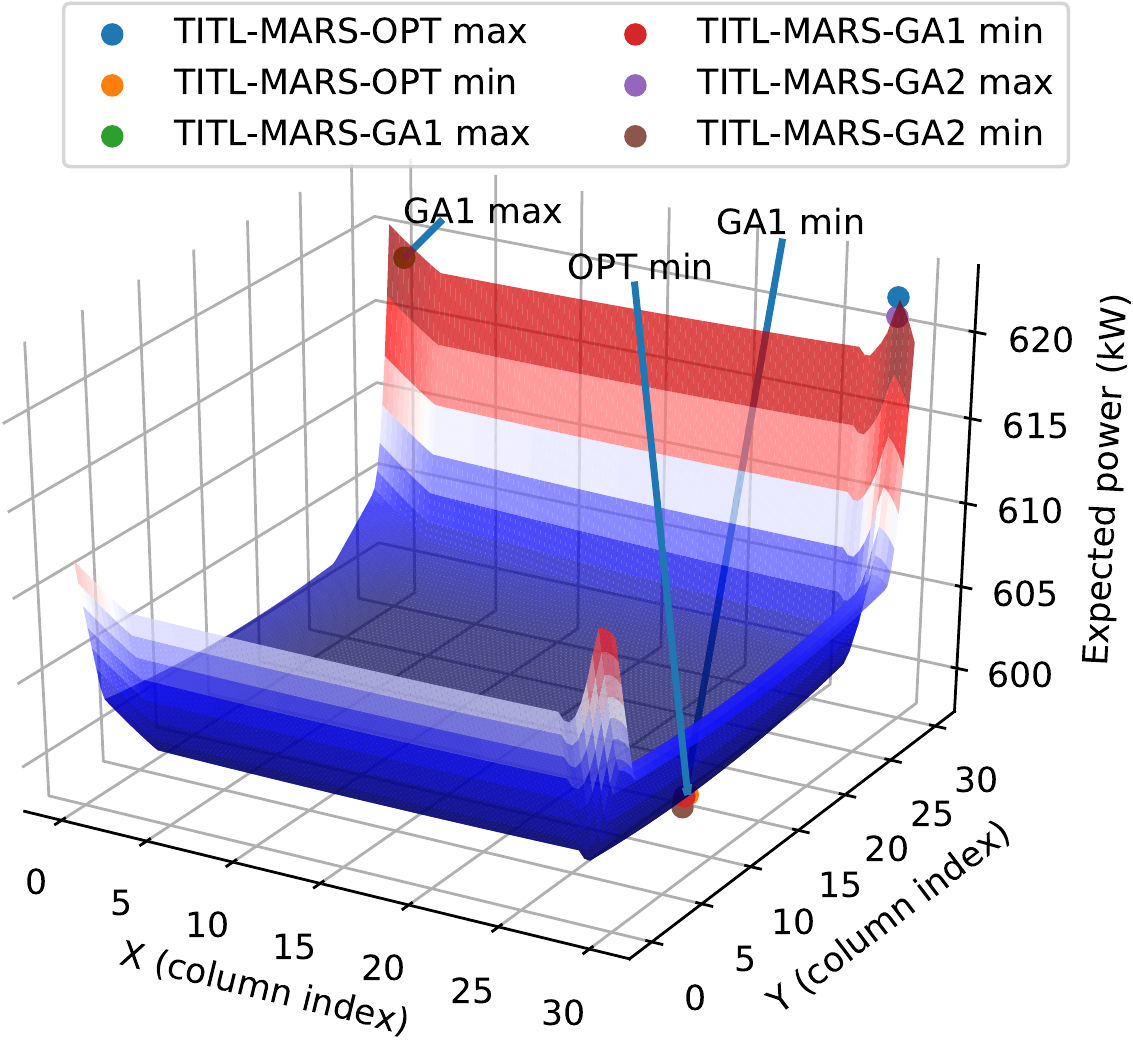}
		\caption{{\small Case 3 : $f_{w3}$}}
		\label{fig:mwd_surface_ex_3}
	\end{subfigure}
	\hspace{0em}
	\begin{subfigure}[b]{0.45\textwidth}
		\centering
		\includegraphics[width=2.2in]{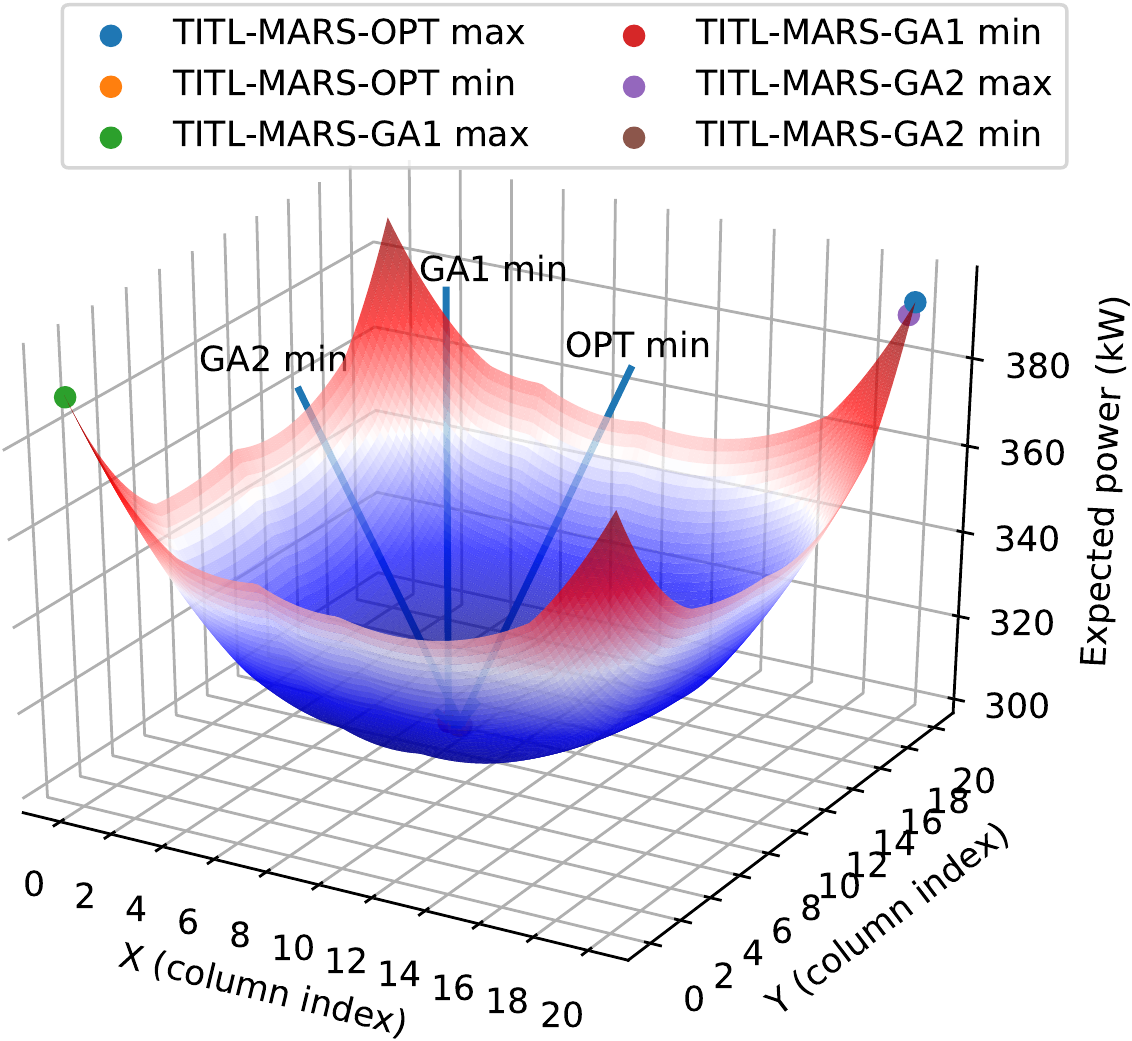}
		\caption[]%
		{{\small Case 4 : $f_{w4}$}}    
		\label{fig:mwd_surface_ex_4}
	\end{subfigure}
    \vspace{0\baselineskip}
	\caption[]
	{\small One run result comparison of TITL-MARS-OPT and TITL-MARS-GA on wind farm power distribution TITL-MARS model} 
	\label{fig:mwd_surfaces}
\end{figure*}

TITL-MARS-OPT and TITL-MARS-GA are used to optimize on the wind farm power distribution models to find a global maximum and a minimum, and the results are shown in Figure~\ref{fig:mwd_surfaces} and summarized in Table~\ref{tab:mars_ti_opt_cmp1}. The results are the average value of 30 executions. The table shows the optimal values derived from the TITL-MARS-OPT and TITL-MARS-GA, as well as the computation time in seconds. The result shows that the TITL-MARS-OPT method finds  better solutions than the genetic algorithm and uses less time. The maximum value and minimum value are very useful before actually building the wind turbines. The maximum location indicates that it is the best location to build a wind turbine based on the given requirements. The minimum location indicates that this location is the worst location on the wind farm, and if the budget is tight, the piece of land around the minimum location can be neglected.

\begin{table}
	\centering
	\small
	\caption{Comparison of TITL-MARS-OPT and TITL-MARS-GA on wind farm power distribution TITL-MARS models}
	\label{tab:mars_ti_opt_cmp1}
	\scalebox{1.0}{
	\begin{tabular}{c|l|cccc}
		\hline
		Function&Measurement&MARS-OPT&MARS-GA 1&MARS-GA 2&MARS-GD\\
		\hline
		\multirow{4}{*}{$f_{w1}$}&\cellcolor{blue!25}Maximum&\cellcolor{blue!25}636.76&\cellcolor{blue!25}618.40&\cellcolor{blue!25}630.75&\cellcolor{blue!25}608.50\\
		&Time(seconds)&	0.15&	1.45&9.06&0.04\\
		\cline{2-4}
		&\cellcolor{green!25}Minimum&\cellcolor{green!25}578.30&\cellcolor{green!25}579.18&\cellcolor{green!25}578.33&\cellcolor{green!25}579.21\\
		&Time(seconds)&	0.32&	1.46&9.08&0.04\\
		\hline
		
		\multirow{4}{*}{$f_{w2}$}&\cellcolor{blue!25}Maximum&\cellcolor{blue!25}588.10&\cellcolor{blue!25}569.96&\cellcolor{blue!25}585.42&\cellcolor{blue!25}582.50\\
		&Time(seconds)&	0.15&	1.46&9.18&0.04\\
		\cline{2-4}
		&\cellcolor{green!25}Minimum&\cellcolor{green!25}545.90&\cellcolor{green!25}545.88&\cellcolor{green!25}545.85&\cellcolor{green!25}545.91\\
		&Time(seconds)&0.40&1.48&9.15&0.08\\
		
		\hline
		\multirow{4}{*}{$f_{w3}$}		&\cellcolor{blue!25}Maximum&\cellcolor{blue!25}622.10&\cellcolor{blue!25}608.13&\cellcolor{blue!25}621.48&\cellcolor{blue!25}609.20\\
		&Time(seconds)&0.15&1.49&9.36&0.08\\
		\cline{2-4}
		&\cellcolor{green!25}Minimum&\cellcolor{green!25}599.50&\cellcolor{green!25}599.77&\cellcolor{green!25}599.51&\cellcolor{green!25}599.86\\
		&Time(seconds)&0.85&1.52&9.34&0.05\\
		\hline

		\multirow{4}{*}{$f_{w4}$}
		&\cellcolor{blue!25}Maximum&\cellcolor{blue!25}393.80&\cellcolor{blue!25}387.92&\cellcolor{blue!25}393.13&\cellcolor{blue!25}382.71\\
		&Time(seconds)&0.04&1.57&9.72&0.23\\
		\cline{2-4}
		&\cellcolor{green!25}Minimum&\cellcolor{green!25}303.21&\cellcolor{green!25}303.46&\cellcolor{green!25}303.21&\cellcolor{green!25}303.91\\
		&Time(seconds)&	0.51&1.58&9.67&0.01\\
		\hline
	\end{tabular}}
\end{table}

\subsection{Optimization of other functions}
\begin{figure*}[]
	\centering
	
	\hspace{0em}
	\begin{subfigure}[b]{0.45\textwidth}
		\centering
		\includegraphics[width=2.2in]{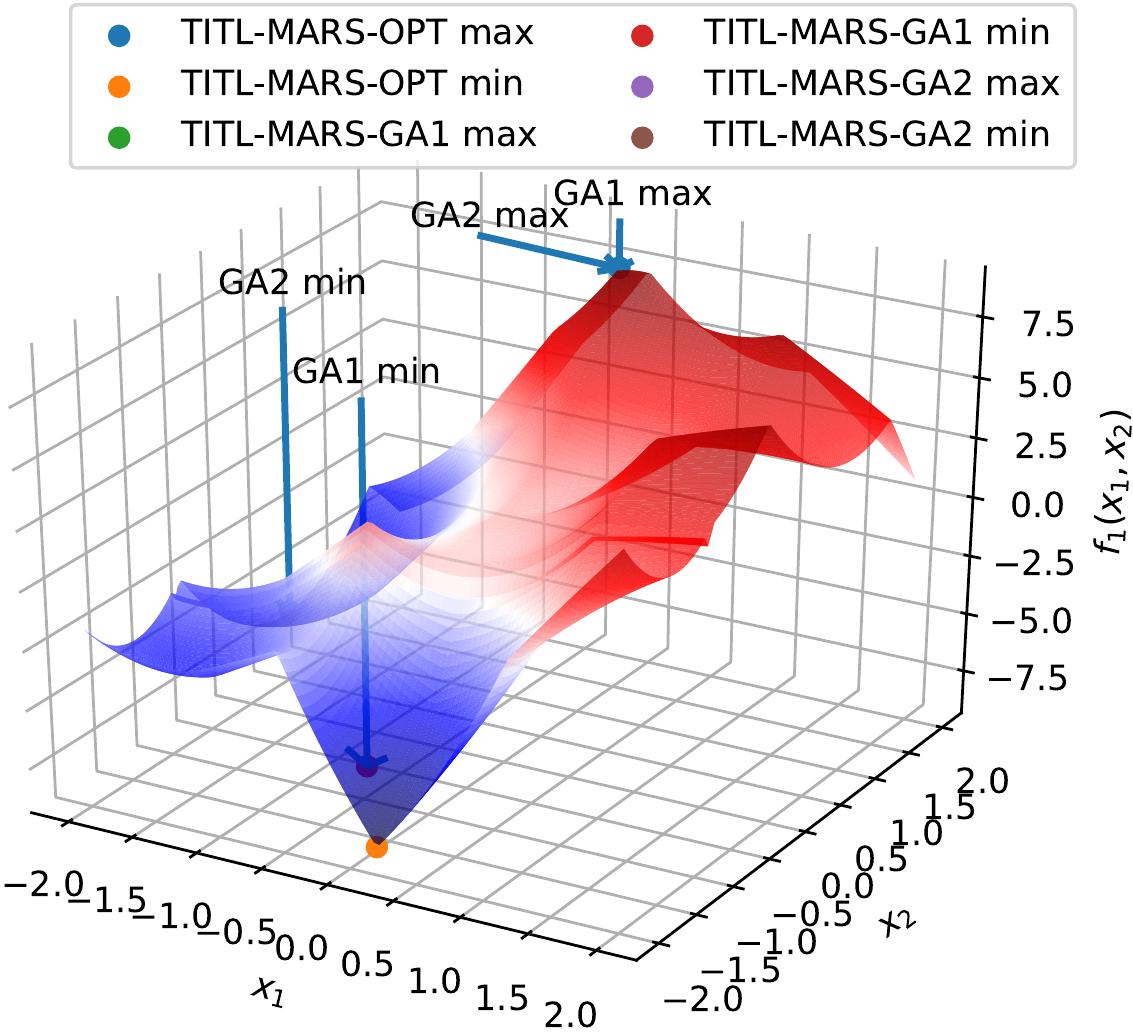}
		\caption[]%
		{{\small $f_1$}}    
		\label{fig:f1_surface_ex}
	\end{subfigure}
	\begin{subfigure}[b]{0.45\textwidth}  
		\centering 
		\includegraphics[width=2.2in]{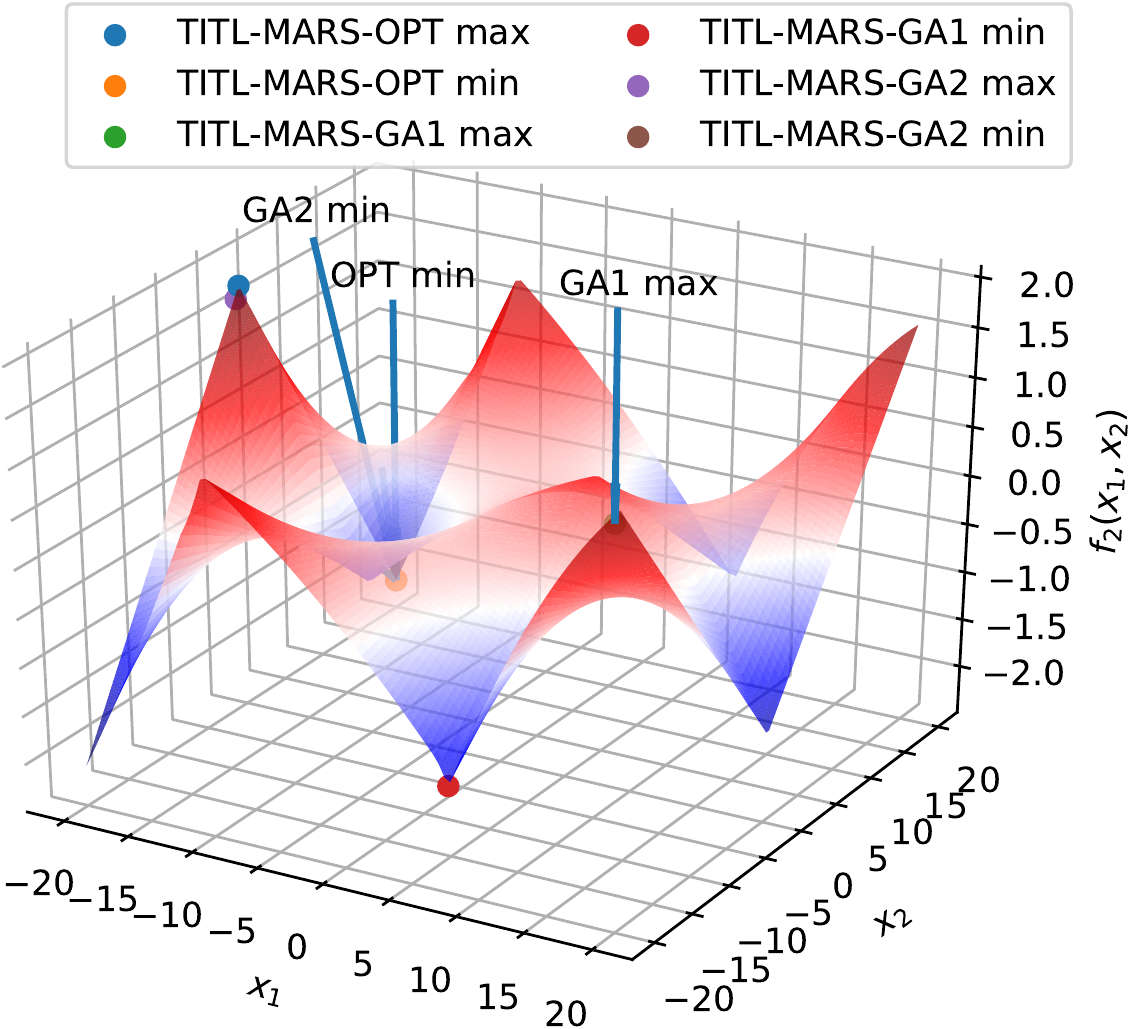}
		\caption[]%
		{{\small $f_2$}}
		\label{fig:f2_surface_ex}
	\end{subfigure}
	
	\caption{\small One run result comparison of TITL-MARS-OPT and TITL-MARS-GA on $f_1$ and $f_2$ TITL-MARS models} 
	\label{fig:other_surfaces}
\end{figure*}

TITL-MARS-OPT method and TITL-MARS-GA are tested to optimize other six TITL-MARS models to find the global maximum and minimum. The first two TITL-MARS models $f_1$ and $f_2$ are two-dimensional~\citep{miyata2005free}.  The $f_3$ and $f_4$ are 10-dimensional TITL-MARS models. $f_5$ is 19-dimensional, and $f_6$ is 21-dimensional~\citep{ariyajunya2013adaptive}. The results are shown in Figure~\ref{fig:other_surfaces} and summarized in Table~\ref{tab:mars_ti_opt_cmp2}. The results are the average value of 30 runs and show TITL-MARS-OPT is superior to TITL-MARS-GA. The result shows that TITL-MARS-OPT achieves better solutions than TITL-MARS-GA and uses less time, which is consistent with the prior result. The results also show that TITL-MARS-OPT is robust in dealing both low-dimensional and high-dimensional TITL-MARS models.
\begin{table}
	\centering
	\small
	\caption{Result comparison of TITL-MARS-OPT and TITL-MARS-GA on six other TITL-MARS mathematical models}
	\label{tab:mars_ti_opt_cmp2}
	\scalebox{1}{
	\begin{tabular}{c|l|cccc}
		\hline
		Function&Measurement&MARS-OPT&MARS-GA 1&MARS-GA 2&MARS-GD\\
		\hline

		\multirow{4}{*}{$f_{1}$}&\cellcolor{blue!25}Maximum&\cellcolor{blue!25}8.30&\cellcolor{blue!25}6.23&\cellcolor{blue!25}8.27&\cellcolor{blue!25}6.60\\
		&Time(seconds)&	0.70&1.50&9.25&2.38\\
		\cline{2-4}
		&\cellcolor{green!25}Minimum&\cellcolor{green!25}-8.20&\cellcolor{green!25}-6.42&\cellcolor{green!25}-6.33&\cellcolor{green!25}-5.06\\
		&Time(seconds)&	1.65&1.51 &9.18&0.57\\
		\hline
		
		\multirow{4}{*}{$f_{2}$}
		&\cellcolor{blue!25}Maximum&\cellcolor{blue!25}1.81&\cellcolor{blue!25}1.02 &\cellcolor{blue!25}1.45&\cellcolor{blue!25}1.24\\
		&Time(seconds)&	0.31&1.42&8.91&1.71\\
		\cline{2-4}
		&\cellcolor{green!25}Minimum&\cellcolor{green!25}-2.20&\cellcolor{green!25}-1.38 &\cellcolor{green!25}-2.20&\cellcolor{green!25}-1.38\\
		&Time(seconds)&0.32&1.42 &8.88&1.54\\
		\hline
		
		\multirow{4}{*}{$f_{3}$}
		&\cellcolor{blue!25}Maximum&\cellcolor{blue!25} 5,774.08&\cellcolor{blue!25}5,092.52 &\cellcolor{blue!25}5,410.16&\cellcolor{blue!25} 5661.76\\
		&Time(seconds)&0.02	 &2.71 &15.74 &0.02\\
		\cline{2-4}
		&\cellcolor{green!25}Minimum&\cellcolor{green!25} -1126.39&\cellcolor{green!25}2,289.47 &\cellcolor{green!25}903.69&\cellcolor{green!25}-806.86 \\
		&Time(seconds)&0.02 &2.67 &15.81&10.35\\
		\hline
		
		\multirow{4}{*}{$f_{4}$}
		&\cellcolor{blue!25}Maximum&\cellcolor{blue!25} 48,800.26&\cellcolor{blue!25}-13,516.73 &\cellcolor{blue!25}-19,762.80&\cellcolor{blue!25}-1,891,678.92\\
		&Time(seconds)&96.32&2.98 &17.67&14.24\\
		\cline{2-4}
		&\cellcolor{green!25}Minimum&\cellcolor{green!25} -3,952,146.24&\cellcolor{green!25}-2,605,904.95&\cellcolor{green!25} -3,599,093.85&\cellcolor{green!25}-3,119,142.40\\
		&Time(seconds)&0.41 &2.95 &17.57 &0.01\\
		\hline
		
		\multirow{4}{*}{$f_{5}$}
		&\cellcolor{blue!25}Maximum&\cellcolor{blue!25} 97,679.99&\cellcolor{blue!25}78,263.13 &\cellcolor{blue!25}92,458.62&\cellcolor{blue!25}80,461.36\\
		&Time(seconds)&0.02	 &5.20 &30.01&28.45 \\
		\cline{2-4}
		&\cellcolor{green!25}Minimum&\cellcolor{green!25} -15,439.62&\cellcolor{green!25}33,298.80 &\cellcolor{green!25} 5,904.35&\cellcolor{green!25}719.07\\
		&Time(seconds)&2.46 &5.21 &30.38& 28.72\\
		\hline
		
		\multirow{4}{*}{$f_{6}$}
		&\cellcolor{blue!25}Maximum&\cellcolor{blue!25} 111,225.22&\cellcolor{blue!25} 63,506.19&\cellcolor{blue!25}105,645.30&\cellcolor{blue!25}88956.05\\
		&Time(seconds)&0.02	 &5.34 &30.62& 29.90\\
		\cline{2-4}
		&\cellcolor{green!25}Minimum&\cellcolor{green!25} -14,215.61&\cellcolor{green!25}17,886.61&\cellcolor{green!25}494.43&\cellcolor{green!25}-5757.07\\
		&Time(seconds)&0.04 &5.33 &30.58&12.64\\
		\hline
		
	\end{tabular}}
\end{table}

\begin{table}
	\centering
	\small
	\caption{Result comparison of true function value}
	\label{tab:mars_true_value_cmp__}
	\scalebox{0.6}{
		\begin{tabular}{c|l|ccccccc}
			\hline
			Function&Measurement&OPT function(model)&MARS+GD&GA 1&GA 2&GD&PL-OPT&PL OPT+GD\\
			\hline

			\multirow{4}{*}{$f_{1}$}&\cellcolor{blue!25}Maximum&\cellcolor{blue!25}8.01(8.32)&\cellcolor{blue!25}8.21 &\cellcolor{blue!25}8.20&\cellcolor{blue!25}8.21&\cellcolor{blue!25}5.76&\cellcolor{blue!25}-4.76(-5.15)&\cellcolor{blue!25}2.96\\
			&Time(seconds)&0.25	&0.24&0.78&7.58&0.01&0.05&0.00\\
			\cline{2-4}
			&\cellcolor{green!25}Minimum&\cellcolor{green!25}-4.56(-8.16)&\cellcolor{green!25}-6.20&\cellcolor{green!25}-3.75&\cellcolor{green!25}-5.76&\cellcolor{green!25}-3.30&\cellcolor{green!25}-4.86(-5.58)&\cellcolor{green!25}-6.20\\
			&Time(seconds)&10.57&0.73	&0.77 &7.54&0.01&0.05&0.00\\
			\hline
			
			\multirow{4}{*}{$f_{2}$}&\cellcolor{blue!25}Maximum&\cellcolor{blue!25} 0.86(1.81)&\cellcolor{blue!25}1.00&\cellcolor{blue!25}1.00&\cellcolor{blue!25}1.00&\cellcolor{blue!25}1.00&\cellcolor{blue!25}0.31(0.56)&\cellcolor{blue!25}1.00\\
			&Time(seconds)&0.10&0.14	&0.59&5.25&0.10&0.24&0.00\\
			\cline{2-4}
			&\cellcolor{green!25}Minimum&\cellcolor{green!25}-0.61(-2.20)&\cellcolor{green!25}-1.00&\cellcolor{green!25}-0.87&\cellcolor{green!25}-1.00&\cellcolor{green!25}-0.96&\cellcolor{green!25}-0.46(-0.53)&\cellcolor{green!25}-1.00\\
			&Time(seconds)&2.65&0.22	& 0.54&5.27&0.10&0.13&0.00\\
			\hline
			
			\multirow{4}{*}{$f_{3}$}&\cellcolor{blue!25}Maximum&\cellcolor{blue!25} 6,029.13(5,774.08)&\cellcolor{blue!25}6132.00&\cellcolor{blue!25}5,488.40&\cellcolor{blue!25}5,647.09&\cellcolor{blue!25}5893.00&\cellcolor{blue!25}4821.00(4667.86)&\cellcolor{blue!25}6132.00\\
			&Time(seconds)&0.01	&0.01&1.90&18.81&0.01&0.07&0.00\\
			\cline{2-4}
			&\cellcolor{green!25}Minimum&\cellcolor{green!25}-826.08(-1,126.39)&\cellcolor{green!25}-923.86&\cellcolor{green!25}2,755.59&\cellcolor{green!25}742.00&\cellcolor{green!25}-525.09&\cellcolor{green!25}288.27(118.11)&\cellcolor{green!25}-123.44\\
			&Time(seconds)&	0.01&0.01& 1.85&18.81&0.47&0.05&0.00\\
			\hline
			
			\multirow{4}{*}{$f_{4}$}&\cellcolor{blue!25}Maximum&\cellcolor{blue!25}-1,991.28(48,800.26)&\cellcolor{blue!25}-0.01 &\cellcolor{blue!25} -1,643,753.96&\cellcolor{blue!25}-131,067.27&\cellcolor{blue!25}-0.01&\cellcolor{blue!25}-505,071.14(-398,359.99)&\cellcolor{blue!25}-0.01\\
			&Time(seconds)&107.47&112.74	&2.12&21.98&2.47&0.07&0.00\\
			\cline{2-4}
			&\cellcolor{green!25}Minimum&\cellcolor{green!25}-4,581,942.36(-3,952,146.24)&\cellcolor{green!25}-4,616,810.03&\cellcolor{green!25}-3,651,963.86&\cellcolor{green!25}-3,834,009.84&\cellcolor{green!25}-3,161,354.60&\cellcolor{green!25}-4,433,291.39(-3,945,225.00)&\cellcolor{green!25}-4,616,810.03\\
			&Time(seconds)&0.40&0.42	&2.04 &21.83&0.01&0.10&0.00\\
			\hline
			
	\end{tabular}}
\end{table}

\begin{figure*}[]
	\centering
	
	\hspace{0em}
	\begin{subfigure}[b]{0.23\textwidth}
		\centering
		\includegraphics[height=1.4in]{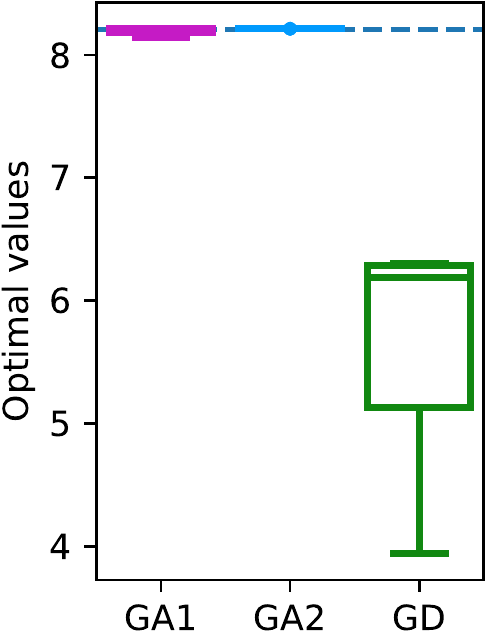}
		\caption[]%
		{{\small $f_1$}}    
		\label{fig:f1_max}
	\end{subfigure}
	\begin{subfigure}[b]{0.23\textwidth}  
		\centering 
		\includegraphics[height=1.4in]{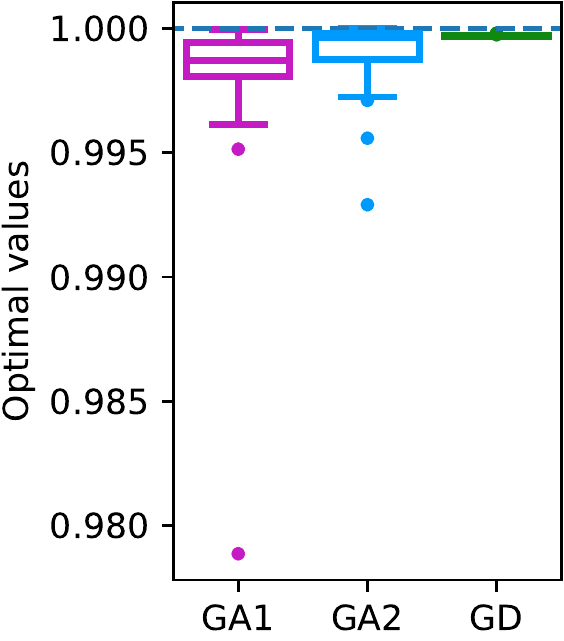}
		\caption[]%
		{{\small $f_2$}}
		\label{fig:f2_max}
	\end{subfigure}
	\begin{subfigure}[b]{0.23\textwidth}
		\centering
		\includegraphics[height=1.4in]{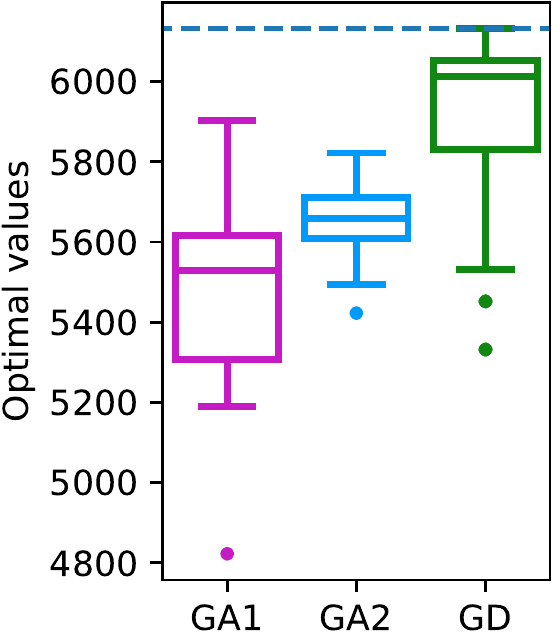}
		\caption[]%
		{{\small $f_3$}}    
		\label{fig:f3_max}
	\end{subfigure}
	\begin{subfigure}[b]{0.23\textwidth}  
		\centering 
		\includegraphics[height=1.4in]{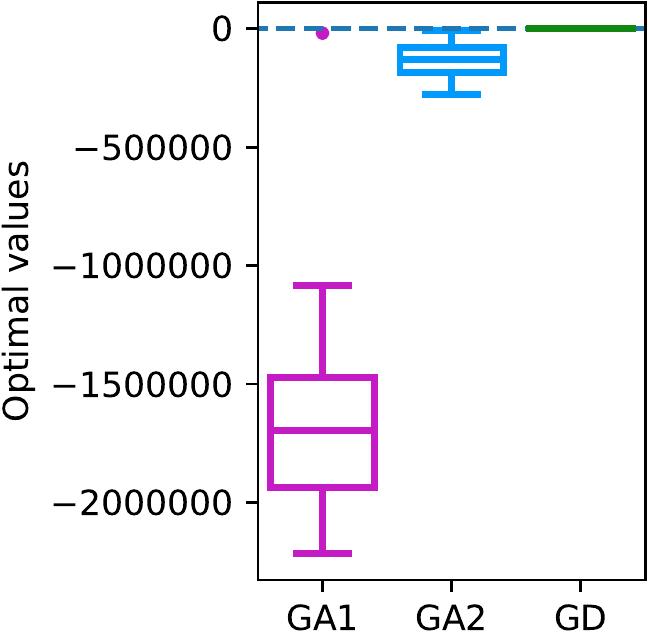}
		\caption[]%
		{{\small $f_4$}}
		\label{fig:f4_max}
	\end{subfigure}
	
	\caption{\small maximum values boxplot} 
	\label{fig:max_bp}
\end{figure*}

\begin{figure*}[]
	\centering
	
	\hspace{0em}
	\begin{subfigure}[b]{0.23\textwidth}
		\centering
		\includegraphics[height=1.4in]{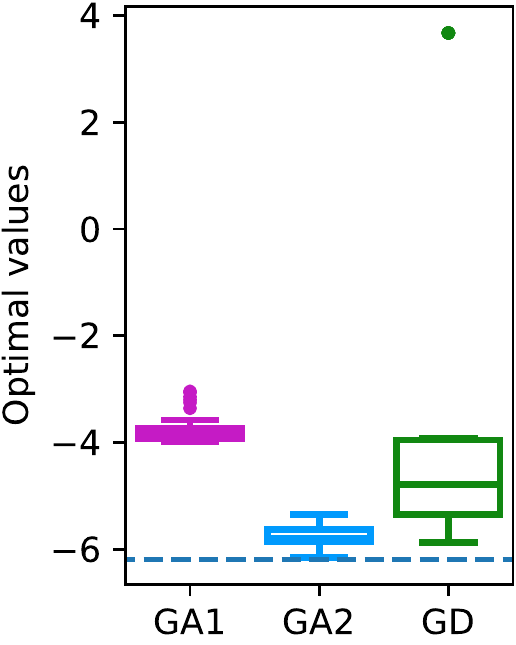}
		\caption[]%
		{{\small $f_1$}}    
		\label{fig:f1_min}
	\end{subfigure}
	\begin{subfigure}[b]{0.23\textwidth}  
		\centering 
		\includegraphics[height=1.4in]{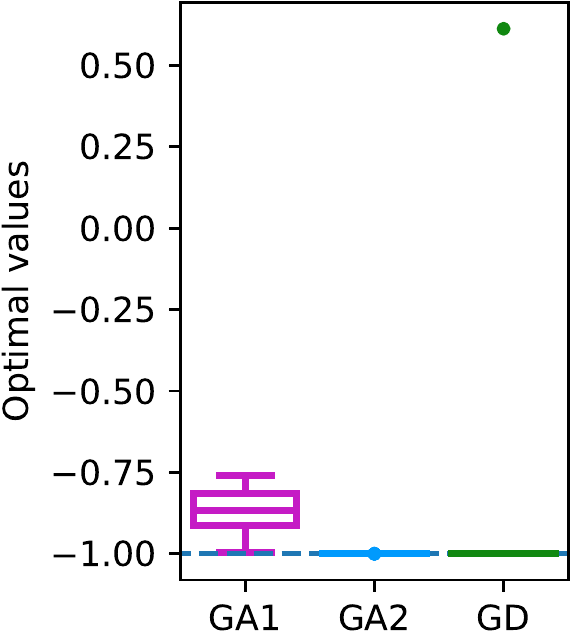}
		\caption[]%
		{{\small $f_2$}}
		\label{fig:f2_min}
	\end{subfigure}
	\begin{subfigure}[b]{0.23\textwidth}
		\centering
		\includegraphics[height=1.4in]{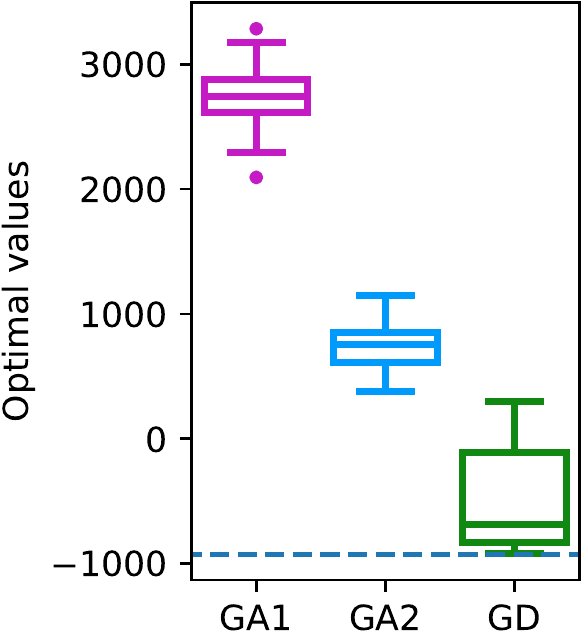}
		\caption[]%
		{{\small $f_3$}}    
		\label{fig:f3_min}
	\end{subfigure}
	\begin{subfigure}[b]{0.23\textwidth}  
		\centering 
		\includegraphics[height=1.4in]{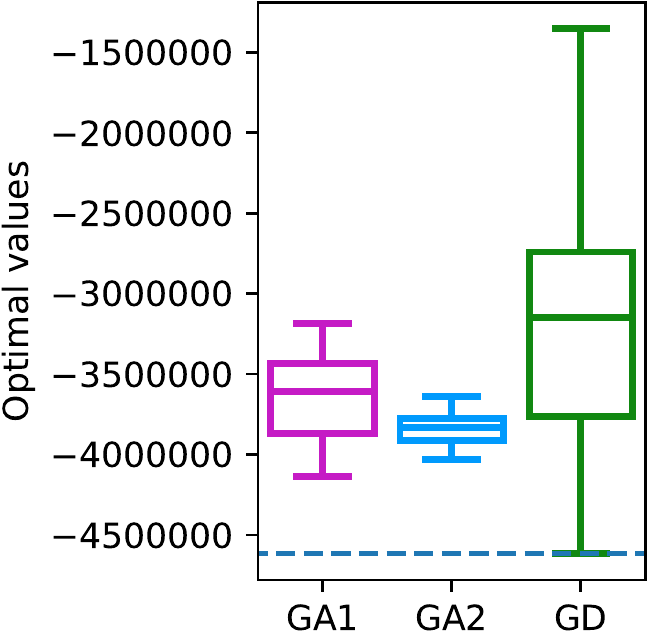}
		\caption[]%
		{{\small $f_4$}}
		\label{fig:f4_min}
	\end{subfigure}
	
	\caption{\small minimum values boxplot} 
	\label{fig:min_bp}
\end{figure*}

\subsection{$f_1$, $f_2$, $f_3$ and $f_4$ functions}

\begin{flalign}
f_1({x_1,x_2})=& 3(1-x_1)^2\exp(-x_1^2-(x_2+1)^2) -10(\frac{x_1}{5}-x_1^3-x_2^5)\exp(-x_1^2-x_2^2)&&\\\nonumber
&-\frac{1}{3}\exp(-(x_1+1)^2-x_2^2)+2x_1, &&\\\nonumber
&-2\leqslant x_1 \leqslant 2,-2\leqslant x_2 \leqslant 2&&\nonumber
\end{flalign}

\begin{flalign}
f_2(x_1,x_2)=&\sin \left(\frac{\pi x_1}{12}\right) \cos \left(\frac{\pi x_2}{16}\right)&&\\&-20\leqslant x_1\leqslant20,-20\leqslant x_2\leqslant 20&&\nonumber
\end{flalign}

\begin{flalign}
f_3(\mathbf{x})=&x_1^2+x_2^2+x_1x_2-14x_1-16x_2+(x_3-10)^2-4(x_4-5)^2+(x_5-3)^2&&\\\nonumber
&+2(x_6-1)^2+5x_7^2+7(x_8-11)^2+2(x_9-10)^2+2(x_{10}-7)^2+45&&\\
&-10\leqslant x_i \leqslant 10&&\nonumber
\end{flalign}

\begin{flalign}
f_4(\mathbf{x})=&\sum\limits_{j=1}^{10}\exp(x_j)\left(c_j+x_j-\ln\sum\limits_{k=1}^{10}\exp(x_k)\right)&&\\\nonumber
&\mathbf{c}=[-0.6089,-17.164,-34.054,-5.914,-24.721,-14.986,-24.100,-10.708,&&\\\nonumber&-26.662,-22.179]&&\\
&-10\leqslant x_i \leqslant 10&&\nonumber
\end{flalign}

\begin{figure*}
	\centering
	\begin{subfigure}[b]{0.45\textwidth}  
		\centering 
		\includegraphics[height=1.7in]{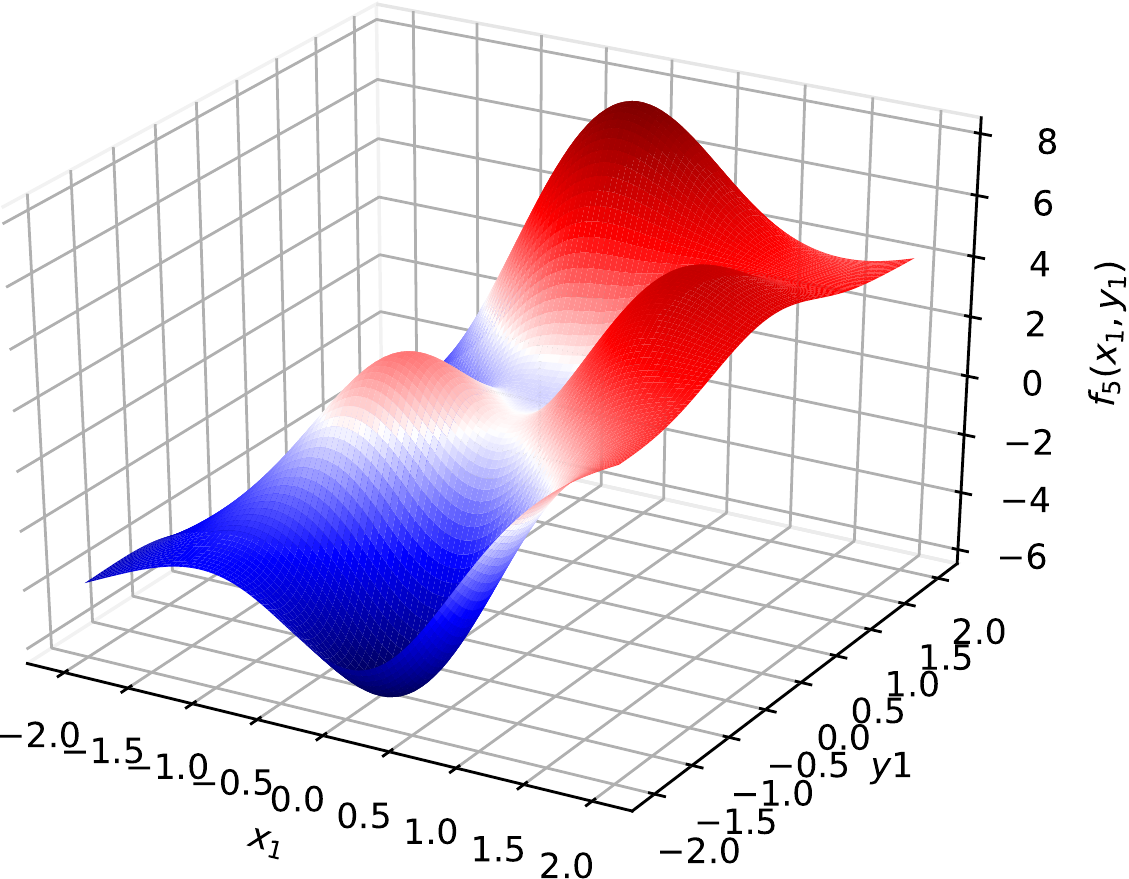}
		\caption[]%
		{{\small $f_1(x_1,x_2)$}}
		\label{fig:f1_s}
	\end{subfigure}
	\hspace{0em}
	\begin{subfigure}[b]{0.45\textwidth}
		\centering
		\includegraphics[height=1.7in]{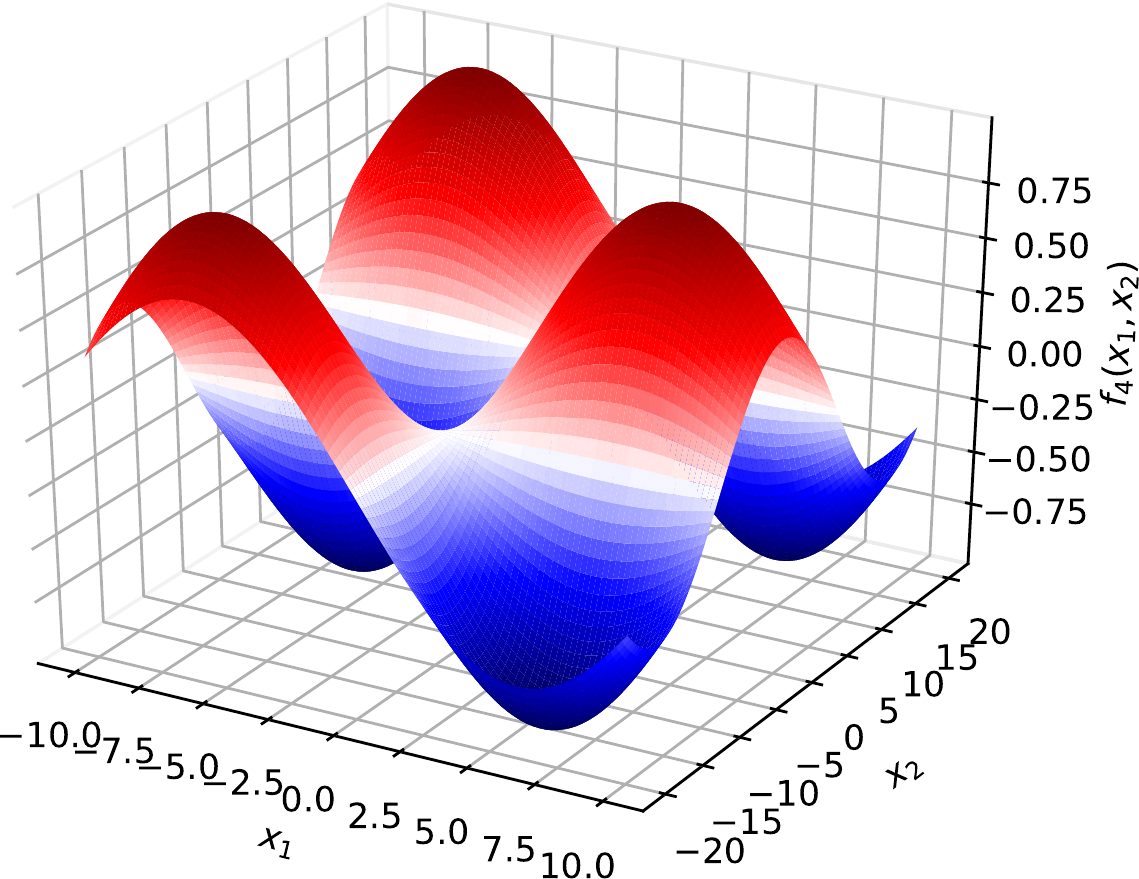}
		\caption[]%
		{{\small $f_2(x_1,x_2)$}}    
		\label{fig:f2_surface}
	\end{subfigure}
	\vspace{0\baselineskip}
	\begin{subfigure}[b]{0.45\textwidth}  
		\centering 
		\includegraphics[height=1.7in]{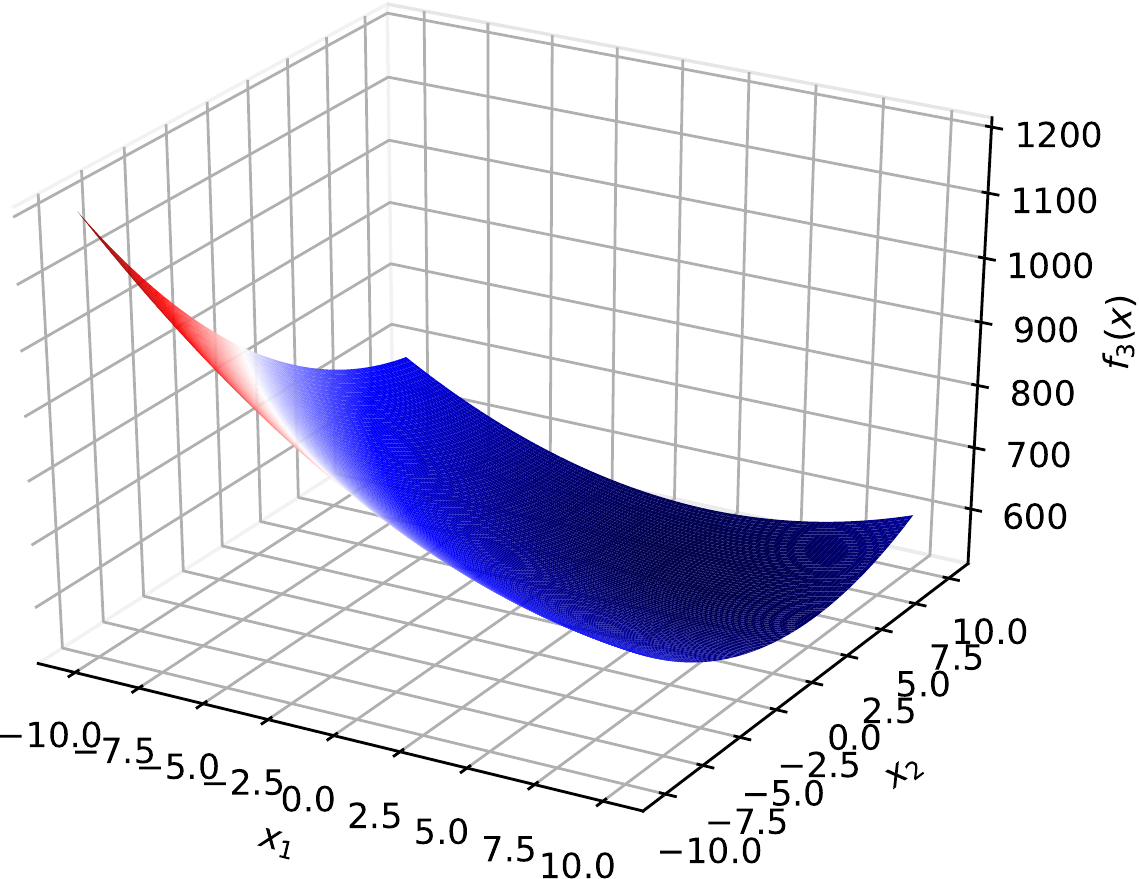}
		\caption[]%
		{{\small $f_3(x_1,x_2) \mbox{ other } x_i=4.0$}}
		\label{fig:f3_s}
	\end{subfigure}
	\hspace{0em}	
	\begin{subfigure}[b]{0.45\textwidth}
		\centering
		\includegraphics[height=1.7in]{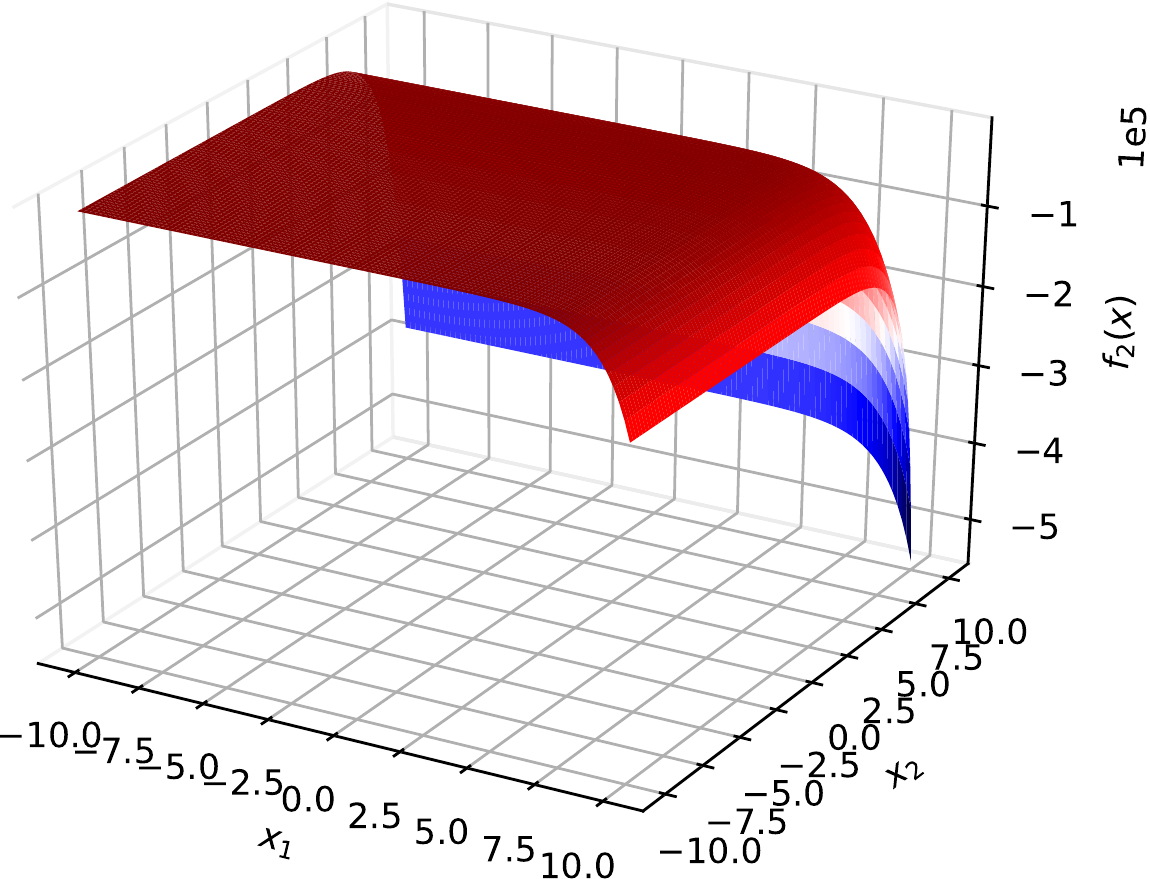}
		\caption[]%
		{{\small $f_4(x_1,x_2) \mbox{ other } x_i=3.0$}}     
		\label{fig:f4_s}
	\end{subfigure}
	\vspace{0\baselineskip}
	\caption{\small Surfaces of dataset functions} 
	\label{fig:f_surfaces}
\end{figure*}

\section{Conclusion}\label{conclusion}
In this paper, a new method (TITL-MARS-OPT) is proposed to globally optimize analytically on the two-way interaction truncated linear MARS (TITL-MARS) by using mixed integer quadratic programming. We verified the presented TITL-MARS-OPT method on the wind farm power distribution TITL-MARS models and six other mathematical TITL-MARS models. The application on wind farm power distribution models gives the best location and worst location information on the wind farm. The testing TITL-MARS models are from 2-dimensions to up to 21-dimensions, and it shows the TITL-MARS-OPT method is robust in dealing with TITL-MARS models with varied dimensions. We also compared the TITL-MARS-OPT method with TITL-MARS-GA in TITL-MARS model optimization, and it shows that the new method can achieve better accuracy and time efficiency. TITL-MARS-OPT can achieve as high as 316\% and on average 46\% better solution quality and is on average 175\% faster than TITL-MARS-GA. In addition, the Python code and the testing models of this paper are made open source, and it will contribute to the study of TITL-MARS models and optimization.

\bibliographystyle{unsrtnat}
\bibliography{main}

\appendix
\section{Supplemental materials}
\subsection{$f_1$, $f_2$, $f_3$ and $f_4$ functions}

\begin{flalign}
f_1({x_1,x_2})=& 3(1-x_1)^2\exp(-x_1^2-(x_2+1)^2) -10(\frac{x_1}{5}-x_1^3-x_2^5)\exp(-x_1^2-x_2^2)&&\\\nonumber
&-\frac{1}{3}\exp(-(x_1+1)^2-x_2^2)+2x_1, &&\\\nonumber
&-2\leqslant x_1 \leqslant 2,-2\leqslant x_2 \leqslant 2&&\nonumber
\end{flalign}

\begin{flalign}
f_2(x_1,x_2)=&\sin \left(\frac{\pi x_1}{12}\right) \cos \left(\frac{\pi x_2}{16}\right)&&\\&-20\leqslant x_1\leqslant20,-20\leqslant x_2\leqslant 20&&\nonumber
\end{flalign}

\begin{flalign}
f_3(\mathbf{x})=&x_1^2+x_2^2+x_1x_2-14x_1-16x_2+(x_3-10)^2-4(x_4-5)^2+(x_5-3)^2&&\\\nonumber
&+2(x_6-1)^2+5x_7^2+7(x_8-11)^2+2(x_9-10)^2+2(x_{10}-7)^2+45&&\\
&-10\leqslant x_i \leqslant 10&&\nonumber
\end{flalign}

\begin{flalign}
f_4(\mathbf{x})=&\sum\limits_{j=1}^{10}\exp(x_j)\left(c_j+x_j-\ln\sum\limits_{k=1}^{10}\exp(x_k)\right)&&\\\nonumber
&\mathbf{c}=[-0.6089,-17.164,-34.054,-5.914,-24.721,-14.986,-24.100,-10.708,&&\\\nonumber&-26.662,-22.179]&&\\
&-10\leqslant x_i \leqslant 10&&\nonumber
\end{flalign}

\subsection{$f_5$ and $f_6$ functions}
The datasets to generate $f_5$ and $f_6$ are from~\cite{ariyajunya2013adaptive}. \cite{ariyajunya2013adaptive} applied adaptive dynamic programming for high-dimensional, multicollinear state sapce and used an Atlanta ozone pollution problem as the case study. The datasets used in this paper are from the fourth stage and the third stage with low variance inflation factors.

\subsection{Details about transforming TITL-MARS into MIQP}
The detailed steps of transforming TITL-MARS into MIQP are given in Algorithm~\ref{alg:MARS_2way_2_miqp}.
The objective function is the general form of MIQP given in~\eqref{eq:miqpp}.
In step 1, let $\mathbf{x}$ be the $V$ dimensional decision variable of the original MARS model $\hat{f}(\mathbf{x})$. The step 4 defines $D$ as the dimension of the decision variable of the new MIQP problem. The step 6 defines $\mathbf{z}$ as the new decision variable of the objective function~\eqref{eq:miqpp} where the first element is 1.0 for $a_0$ in the original problem.
In step 7, let $\mathbf{c}$ be the coefficient vector for the linear elements in the new MIQP problem which is initialized to be a $\mathbf{0}$ vector. Steps 8 to 11 determine the values of $\mathbf{c}$. The coefficient for the first decision element 1.0 in $\mathbf{z}$ is $a_0$. For the univariate basis function $B_m(\mathbf{x})$, the corresponding element in the decision variable of MIQP is $\eta_{1,m}$ and the coefficient is $a_m$.
Steps 12 to 17 define 
$\mathbf{Q}$ as the $D$ by $D$ coefficient matrix for the quadratic element in MIQP which is symmetric and is initialized to be $\mathbf{0}$. For the two-way interaction basis function $B_m(\mathbf{x})$, the corresponding coefficient in $Q$ matrix is $a_m$. Steps 18 to 31 transform the quadratic terms into linear constraints. 
\begin{algorithm}
	\small
	\caption{Formulation of two-way interaction truncated linear MARS into mixed integer quadratic programming}\label{alg:MARS_2way_2_miqp}
	\KwData{$\hat{f}( \mathbf{x})=a_0+\sum_{m=1}^{M}\left\{a_m \cdot \prod_{k=1}^{Km}[s_{k,m} \cdot (x_{v(k,m)}-t_{v(k,m)}) ]_+ \right\}, K_m \leqslant 2, \,\mathcal{M}$}
	\KwResult{$\hat{f}(\mathbf{z})=\frac{1}{2} \textbf{z}^T\textbf{Q}\textbf{z}+\textbf{c}^T\textbf{z}$ }
	$\textbf{x}=(x_1,x_2,\dots x_v, \dots x_V),\textbf{x} \in \mathbb{R}^V$\\
	$\eta_{k,m}=[s_{k,m} \cdot (x_{v(k,m)}-t_{v(k,m)}) ]_+ ,\, \eta_{k,m} \in \mathbb{R}, \, \eta_{k,m} \geqslant 0$ \\
	$y_{k,m} = \left\{ {\begin{array}{*{10}{ll}}
		{1} & s_{k,m} \cdot (x_{v(k,m)}-t_{v(k,m)}) \geqslant 0\\
		{0}  &  s_{k,m} \cdot (x_{v(k,m)}-t_{v(k,m)})<0
		\end{array}} \right. ,y_{k,m} \in \mathbb{B}$\\
	$D=1+V+\sum_{m=1}^{M-1}K_m$\\
	$\textbf{z} \in \mathbb{R}^D$\\
	$\textbf{z}=(1,x_1, \dots,x_V,\eta_{1,1},y_{1,1},\dots,\underbrace{ \eta_{1,m-1},y_{1,m-1}}_{K_{m-1}=1},\underbrace{ \eta_{1,m},y_{1,m},\eta_{2,m},y_{2,m}}_{K_m=2},\dots,\eta_{K_{{M-1},M-1}},y_{K_{{M-1},M-1}})$\\
	$\textbf{c}=\textbf{0} \in \mathbb{R}^D$\\
	$\textbf{c}_0=a_0$\\
	\For{$m=1\, \mbox{to} \, M$}{\lIf{$K_m=1$}{$\textbf{c}(\eta_{1,m})=a_m$}}
	$\textbf{Q}=\textbf{0} \in \mathbb{R}_{D \times D}$\\
	\For{$m=1\, \mbox{to} \, M$}{\If{$K_m=2$}{$\textbf{Q}(\eta_{1,m},\eta_{2,m})=a_m, \,\, \textbf{Q}(\eta_{2,m},\eta_{1,m})=a_m$}}
	\For{$m=1\, \mbox{to} \, M$}{\If{$K_m=2 \, \mbox{and} \, s_{k,m}= +1$}{Add constraint: $x_{v(k,m)}-\eta_{k,m}-\mathcal{M}\cdot y_{k,m}\geqslant t_{v(k,m)}-\mathcal{M}$\\
			Add constraint: $-x_{v(k,m)}+\eta_{k,m}\geqslant -t_{v(k,m)}$\\
			Add constraint: $-\eta_{k,m}+\mathcal{M}\cdot y_{k,m}\geqslant 0$
	}}
	\For{$m=1\, \mbox{to} \, M$}{\If{$K_m=2 \, \mbox{and} \, s_{k,m}= -1$}{Add constraint: $-x_{v(k,m)}-\eta_{k,m}-\mathcal{M}\cdot y_{k,m}\geqslant -t_{v(k,m)}-\mathcal{M}$\\
			Add constraint: $x_{v(k,m)}+\eta_{k,m}\geqslant t_{v(k,m)}$\\
			Add constraint: $-\eta_{k,m}+\mathcal{M}\cdot y_{k,m}\geqslant 0$
	}}
\end{algorithm}

\end{document}